\documentclass[preprint,11pt,compress]{elsarticle}

\usepackage{amsmath,amssymb,amsthm}
\usepackage{algorithm}
\usepackage{algpseudocode}
\usepackage{booktabs}
\usepackage[table]{xcolor}
\usepackage{hyperref}
\usepackage{cleveref}
\usepackage{lmodern}
\usepackage{microtype}
\usepackage{orcidlink}
\usepackage{ragged2e}
\usepackage{subcaption}
\usepackage{comment}
\usepackage{longtable}

\hypersetup{hidelinks}

\newtheorem{theorem}{Theorem}[section]
\newtheorem{assumption}{Assumption}[section]
\newtheorem{corollary}{Corollary}[section]
\newtheorem{lemma}{Lemma}[section]
\newtheorem{proposition}[theorem]{Proposition}
\theoremstyle{definition}

\theoremstyle{remark}
\newtheorem{remark}{Remark}

\makeatletter
\def\subsection{\@startsection{subsection}{2}{\z@}%
	{-3.25ex \@plus -1ex \@minus -.2ex}%
	{1.5ex \@plus .2ex}%
	{\normalfont\normalsize\bfseries}}

\def\subsubsection{\@startsection{subsubsection}{3}{\z@}%
	{-3.25ex \@plus -1ex \@minus -.2ex}%
	{1.5ex \@plus .2ex}%
	{\normalfont\normalsize\bfseries}}
\makeatother
\newcommand{\R}{\mathbb{R}}

\newcommand{\x}{\mathbf{x}}
\newcommand{\g}{\mathbf{g}}
\newcommand{\dvec}{\mathbf{d}}
\newcommand{\svec}{\mathbf{s}}
\newcommand{\yvec}{\mathbf{y}}

\newcommand{\bu}{\mathbf{u}}
\newcommand{\bv}{\mathbf{v}}
\newcommand{\bz}{\mathbf{z}}
\newcommand{\bh}{\mathbf{h}}
\newcommand{\br}{\mathbf{r}}
\newcommand{\ba}{\mathbf{a}}
\newcommand{\bb}{\mathbf{b}}

\newcommand{\zvec}{\mathbf{z}}

\newcommand{\bW}{\mathbf{W}}
\newcommand{\bH}{\mathbf{H}}
\newcommand{\bR}{\mathbf{R}}

\newcommand{\bmu}{\boldsymbol{\mu}}
\newcommand{\bsigma}{\boldsymbol{\sigma}}

\newcommand{\norm}[1]{\left\lVert #1 \right\rVert}

\usepackage{tikz}
\usetikzlibrary{arrows.meta,positioning,shapes.geometric,fit}

\begin{document}
	
	\begin{frontmatter}
		
		\title{Safeguarded Learning-Assisted Nonlinear Conjugate Gradient Methods}
		
		\author[lab1]{Morteza Kimiaei\corref{cor1}\,\orcidlink{0000-0002-7973-3770}}
		\ead{morteza.kimiaei@univie.ac.at}
		
		\author[lab2]{Saman Babaie--Kafaki\,\orcidlink{0000-0003-0122-8384}}
		\ead{saman.babaiekafaki@unibz.it}

		\affiliation[lab1]{
			organization={Faculty of Mathematics, University of Vienna},
			addressline={Oskar--Morgenstern--Platz 1, A--1090 Vienna},
			country={Austria}}
		
		\affiliation[lab2]{
			organization={Faculty of Engineering, Free University of Bozen--Bolzano},
			addressline={NOI Techpark, via Bruno Buozzi 1, 39100 Bolzano},
			country={Italy}}

		\cortext[cor1]{Corresponding author}

\begin{abstract}
A learning-assisted framework for \textbf{Nonlinear Conjugate Gradient}
(\texttt{NCG}) methods is developed in which lightweight online prediction is
used to support, rather than replace, the classical optimization process. The
proposed approach retains the standard two-term \texttt{NCG} structure and
introduces two independent online extreme learning machine models. The first
model predicts an \texttt{NCG} coefficient during restart recovery or when
persistent stagnation is detected, whereas the second model is used only
after a line-search failure in finite precision arithmetic to propose a bounded correction to the
deterministic fallback step size. Neural predictions are accepted only when
they satisfy deterministic numerical safeguards, including descent, scale,
and out-of-distribution requirements; otherwise, the corresponding
deterministic recovery or fallback mechanism is retained.  Both learning models are updated online using quantities
already generated along the optimization trajectory and therefore require no
additional function or gradient evaluations for training. These learning
frameworks are applied to six well-known classical \texttt{NCG} methods,
producing corresponding low-storage learning-assisted variants with
deterministic safeguards. Numerical experiments on 480 \texttt{CUTEst}
test problems show improved robustness for all six variants at broadly comparable
evaluation costs. A PDE-constrained conductivity-identification example
further illustrates the potential of the proposed approach to reduce
evaluation costs and computational time for problems with expensive function
and gradient evaluations. These findings suggest that restricted online
learning can improve the robustness and computational performance of
classical \texttt{NCG} methods while preserving their low-storage structure
and deterministic optimization safeguards, and is particularly promising
when function or gradient evaluations are expensive.
\end{abstract}


	\end{frontmatter}
	
%	\linenumbers
	
%\tableofcontents

\begin{sloppypar}

\section{Introduction}
\label{sec:introduction}

An \textbf{Extreme Learning Machine} (\texttt{ELM}) is a learning framework that, in its standard form, uses a feedforward neural network ({\tt NN}) with a single hidden layer whose input weights and biases are randomly assigned and then kept fixed, while the output weights are determined by solving a linear least-squares problem~\cite{Huang2006ELM}. This construction avoids iterative backpropagation through all network parameters, making compact \texttt{ELM} models attractive when frequent fitting and prediction are required. Regularized formulations balance fitting accuracy against the magnitude of the output weights~\cite{Huang2012ELM}, while online sequential extensions allow the predictor to incorporate observations as they become available~\cite{Liang2006OSELM}. Beyond regression and classification, \texttt{ELM} models have been incorporated into optimization algorithms to guide search decisions when exact evaluations are expensive. Pavelski et al.~\cite{Pavelski2016ELMOEAD} combined \texttt{ELM} surrogate models with a decomposition-based multiobjective evolutionary algorithm and differential evolution, obtaining promising results under small evaluation budgets on benchmark problems and an airfoil-design problem. Harada~\cite{Harada2023ELDR} combined \texttt{ELM} with pairwise ranking to estimate the relative quality of candidate solutions, improving the performance of a surrogate-assisted evolutionary algorithm on many tested problems, particularly on higher-dimensional benchmarks. These studies illustrate how inexpensive learned models can support optimization by approximating costly responses or prioritizing promising candidates. They also motivate investigating a complementary role for \texttt{ELM}s in classical numerical methods: learning from the optimization trajectory to assist selected algorithmic decisions.

Consider the unconstrained optimization problem
\begin{equation}
	\label{eq:problem}
	\min_{\x\in\R^n} f(\x),
\end{equation}
which frequently arises in practical applications, where $f:\R^n\rightarrow\R$ is often assumed to be continuously differentiable. Starting from
\(\x_0\in\R^n\), a line-search method generates a sequence of approximations according to
\begin{equation}
	\label{eq:iteration}
	\x_{k+1}:=\x_k+\alpha_k\dvec_k,\qquad k = 0,1,\dots,
\end{equation}
where $\alpha_k>0$ denotes the step length taken along the search direction $\dvec_k$
\cite{Andriebook}. For large-scale cases of \eqref{eq:problem}, \textbf{Nonlinear Conjugate Gradient}
(\texttt{NCG}) methods are particularly attractive because they require the
storage of only a few vectors and avoid storing or factorizing a Hessian
approximation \cite{hager2006survey}. The two-term search direction of an \texttt{NCG}  is formulated by
\begin{equation}
	\label{eq:two-term-direction}
	\dvec_0:=-\g_0, \qquad
	\dvec_k:=-\g_k+\beta_k\dvec_{k-1},\qquad k = 1,2,\dots,
\end{equation}
where \(\g_k:=\nabla f(\x_k)\), and the scalar $\beta_k$, here called the \textbf{Conjugate Gradient} (\texttt{CG})
\textbf{coefficient}, implicitly incorporates information from the previous solution
trajectory through the so-called conjugacy property \cite{SBKRairo}.
%Throughout, \(f_k:=f(\x_k)\) and \(\g(\x):=\nabla f(\x)\), and unadorned vector norms are Euclidean, while the subscript infinity denotes the maximum absolute component.

The apparent simplicity of \eqref{eq:two-term-direction} conceals a challenging design problem. A value of \(\beta_k\) that is too small repeatedly drives the method toward \textbf{Steepest Descent} (\texttt{SD}) and can intensify oscillations across narrow valleys \cite{nocedal1999numerical}. Meanwhile, an excessively large value for \(\beta_k\) can cause a wide angular distance from the \texttt{SD} direction
\begin{equation}
\label{eq:steepest-descent-direction}
\dvec_k=-\g_k,
\end{equation}
and consequently may increase the line-search cost or even violate the descent condition. Moreover, several well-known formulas involve curvature denominators that may become small or unreliable in finite-precision arithmetic \cite{Kimiaei2026}. In such situations, standard restart strategies often safeguard the algorithm by setting \(\beta_k=0\) \cite{DaiLiaoLi2004Restart}, but frequent restarts discard useful historical information and may diminish the practical advantages of \texttt{NCG} methods \cite{Kimiaei2026}.

In this study, we consider the Hestenes--Stiefel (\texttt{HS}), Polak--Ribi\`ere--Polyak (\texttt{PRP}), Hager--Zhang (\texttt{HZ}), Dai--Yuan (\texttt{DY}), Fletcher--Reeves (\texttt{FR}), and Dai--Liao (\texttt{DL}) formulas for the \texttt{CG} coefficient \cite{hager2006survey,SBKRairo}, to develop their corresponding \textbf{Neural Network} (\texttt{NN})-based variants, denoted by \texttt{HSNN}, \texttt{PRPNN}, \texttt{HZNN}, \texttt{DYNN}, \texttt{FRNN}, and \texttt{DLNN}, respectively, where \texttt{NN} indicates assistance provided by online \texttt{ELM} models. We keep the original coefficient for each method as the primary
setting in our approach. In particular, whenever the resulting analytical direction satisfies
\begin{equation}
\label{eq:ls-m3-direction-bounds}
\g_k^{\mathsf T}\dvec_k\leq-\delta\norm{\g_k}^2
\quad \text{and} \quad
\norm{\dvec_k}\leq\kappa_d\norm{\g_k},
\end{equation}
with the constants $0<\delta<1$ and $\kappa_d\geq1$, the first inequality in \eqref{eq:ls-m3-direction-bounds} is referred to as the \textbf{Sufficient Descent Condition} (\texttt{SDC}), while the second is referred to as the \textbf{Direction-Scale Condition} (\texttt{DSC}). When the analytical coefficient is finite and nonzero, no method-specific restart condition is active, and the corresponding analytical direction satisfies both the \texttt{SDC} and the \texttt{DSC}, that direction is retained exactly, except when the separate persistent-stagnation mechanism becomes eligible. Otherwise, the analytical coefficient is deterministically safeguarded before the direction is supplied to the \texttt{cvsrch} solver of Mor{\'e} and Thuente \cite{cvsrch}. Note that, to preserve their original characteristics, the traditional \texttt{HS}, \texttt{PRP}, \texttt{HZ}, \texttt{DY}, \texttt{FR}, and \texttt{DL} implementations are not modified by these additional safeguards.

The proposed framework retains the classical \texttt{CG} coefficient as the primary choice and restricts coefficient learning to restart recovery and persistent stagnation. Restart recovery is triggered by a nonfinite or zero coefficient, a nonfinite direction, or a method-specific restart condition, including the safeguards retained for \texttt{DYNN} and \texttt{FRNN}. It seeks a safe nonzero coefficient using an online \texttt{ELM}~\cite{SunELM2014,Wang2022ELM}, recent valid coefficients, or deterministic safe initialization, reverting to the \texttt{SD} direction if recovery fails. Finite analytical coefficients whose directions violate the \texttt{SDC} or \texttt{DSC} conditions are instead projected deterministically onto the admissible set. An initially admissible classical direction is reconsidered only after persistent stagnation in both the objective value and the Euclidean gradient norm. The tentative \texttt{ELM} direction must then satisfy the local recovery and direction-change tests before the common final \texttt{SDC}/\texttt{DSC} safeguard is applied.

The coefficient-\texttt{ELM} learns from valid classical coefficients using twelve dimension-independent scalar features: gradient ratios, normalized inner products, curvature quantities, the current analytical and previously accepted coefficients, and a method identifier. A second, independent \texttt{ELM} acts only after failure of the unchanged \texttt{cvsrch} line search~\cite{cvsrch}. When no returned point can be certified for direct reuse, it may propose a bounded multiplicative enlargement of the original minimal fallback step. This model learns from accepted points returned by \texttt{cvsrch}, including certified cached Armijo points, and every fallback point must pass the original verification before acceptance.

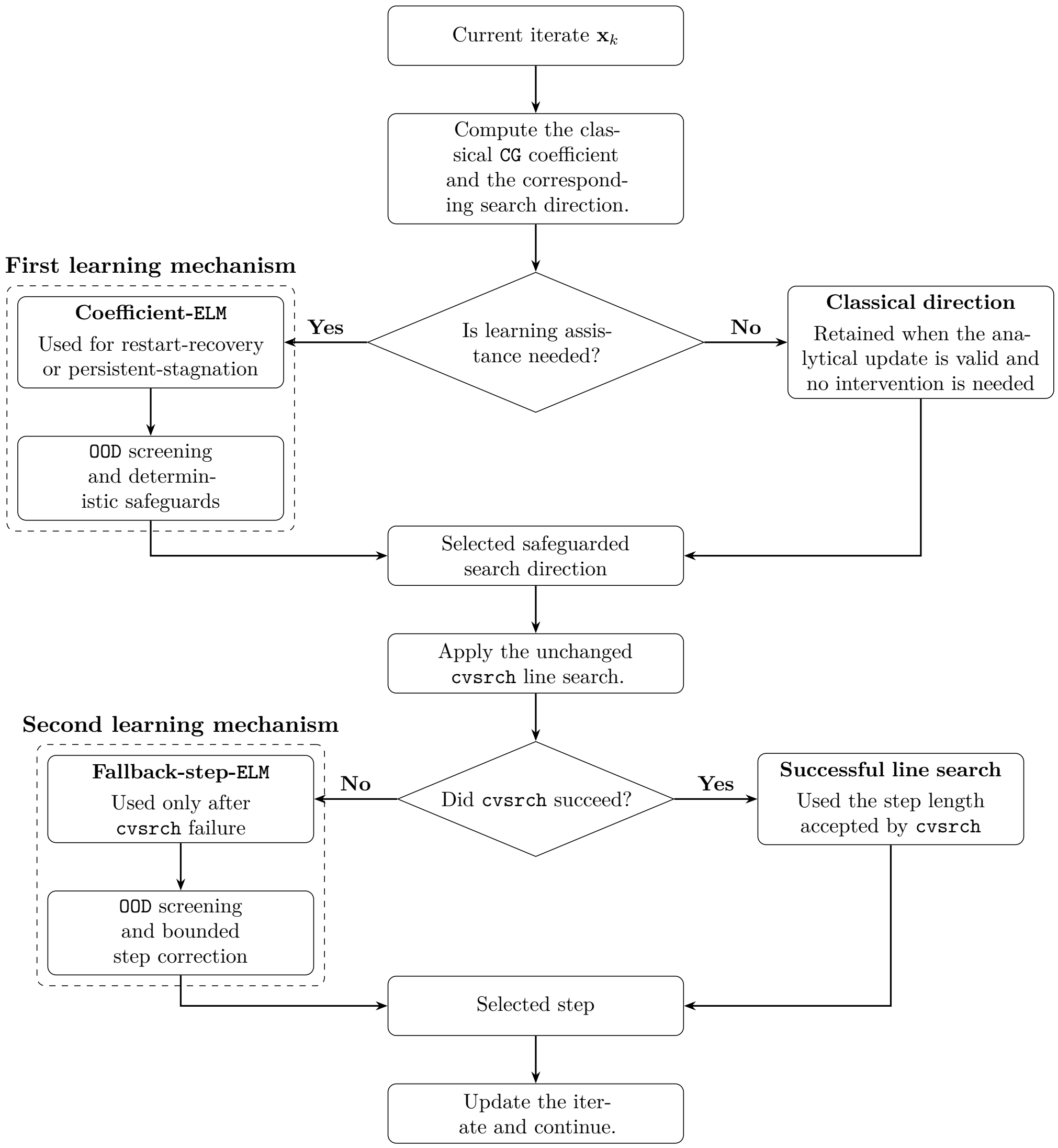
\begin{figure}[htbp]
\centering
\scalebox{0.75}{\begin{tikzpicture}[
	font=\small,
	node distance=8mm and 14mm,
	>=Stealth,
	arrow/.style={
		-{Stealth[length=2mm]},
		thick
	},
	main/.style={
		rectangle,
		draw,
		rounded corners,
		align=center,
		text width=4.7cm,
		minimum height=10mm,
		inner sep=4pt
	},
	side/.style={
		rectangle,
		draw,
		rounded corners,
		align=center,
		text width=4.2cm,
		minimum height=10mm,
		inner sep=4pt
	},
	decision/.style={
		diamond,
		draw,
		align=center,
		aspect=2.4,
		text width=3.6cm,
		inner sep=2pt
	}
]

% =========================================================
% Main CG path
% =========================================================

\node[main] (start)
{Current iterate \(\x_k\)};

\node[main, below=of start] (cg)
{Compute the classical \texttt{CG} coefficient\\
and the corresponding search direction.};

\node[decision, below=of cg] (cgstate)
{Is learning assistance needed?};

\node[main, below=19mm of cgstate] (direction)
{Selected safeguarded search direction};

\node[main, below=of direction] (ls)
{Apply the unchanged \texttt{cvsrch} line search.};

\node[decision, below=of ls] (lsstate)
{Did \texttt{cvsrch} succeed?};

\node[main, below=20mm of lsstate] (step)
{Selected step};

\node[main, below=of step] (update)
{Update the iterate and continue.};

% =========================================================
% First ELM: coefficient
% =========================================================

\node[side, left=14mm of cgstate] (betaelm)
{\textbf{Coefficient-\texttt{ELM}}\\[1mm]
Used for restart-recovery\\
or persistent-stagnation};

\node[side, below=of betaelm] (betacheck)
{\texttt{OOD} screening\\
and deterministic safeguards};

% =========================================================
% Classical branch
% =========================================================

\node[side, right=14mm of cgstate] (classical)
{\textbf{Classical direction}\\[1mm]
Retained when the analytical update
is valid and no intervention is needed};

% =========================================================
% Second ELM: fallback step
% =========================================================

\node[side, left=14mm of lsstate] (stepelm)
{\textbf{Fallback-step-\texttt{ELM}}\\[1mm]
Used only after
\texttt{cvsrch} failure};

\node[side, below=of stepelm] (stepcheck)
{\texttt{OOD} screening\\
and bounded step correction};

% =========================================================
% Successful line search
% =========================================================

\node[side, right=14mm of lsstate] (accepted)
{\textbf{Successful line search}\\[1mm]
Used the step length accepted by
\texttt{cvsrch}};

% =========================================================
% Main arrows
% =========================================================

\draw[arrow] (start) -- (cg);
\draw[arrow] (cg) -- (cgstate);

\draw[arrow]
(cgstate.west)
-- node[above] {\textbf{Yes}}
(betaelm.east);

\draw[arrow]
(betaelm)
-- (betacheck);

\draw[arrow]
(betacheck.south)
|- (direction.west);

\draw[arrow]
(cgstate.east)
-- node[above] {\textbf{No}}
(classical.west);

\draw[arrow]
(classical.south)
|- (direction.east);

\draw[arrow] (direction) -- (ls);
\draw[arrow] (ls) -- (lsstate);

\draw[arrow]
(lsstate.west)
-- node[above] {\textbf{No}}
(stepelm.east);

\draw[arrow]
(stepelm)
-- (stepcheck);

\draw[arrow]
(stepcheck.south)
|- (step.west);

\draw[arrow]
(lsstate.east)
-- node[above] {\textbf{Yes}}
(accepted.west);

\draw[arrow]
(accepted.south)
|- (step.east);

\draw[arrow] (step) -- (update);

% =========================================================
% Labels
% =========================================================

\node[
	draw,
	dashed,
	rounded corners,
	fit=(betaelm)(betacheck),
	inner sep=5pt,
	label={[font=\bfseries]above:First learning mechanism}
] {};

\node[
	draw,
	dashed,
	rounded corners,
	fit=(stepelm)(stepcheck),
	inner sep=5pt,
	label={[font=\bfseries]above:Second learning mechanism}
] {};

\end{tikzpicture}}

\caption{Bird's-eye view of the learning-assisted nonlinear
\texttt{CG} framework.  The classical \texttt{CG} method remains the main
optimization procedure.  The coefficient-\texttt{ELM} assists the construction
of the search direction only during restart recovery or persistent stagnation,
whereas the independent fallback-step-\texttt{ELM} is considered only after
failure of \texttt{cvsrch}.  Both learning branches return to the main
optimization path after \texttt{OOD} screening and the corresponding
deterministic safeguards.}
\label{fig:elm-cg-birdseye}
\end{figure}

Both models are trained online from existing iteration data, without offline candidate searches, external teacher coefficients, or additional function or gradient evaluations for training. Predictions undergo out-of-distribution (\texttt{OOD}) screening against their respective training archives and must satisfy the corresponding deterministic safeguards; unavailable or inadmissible predictions leave the classical direction or applicable deterministic recovery/fallback in control. Figure~\ref{fig:elm-cg-birdseye} summarizes the interaction between these components.

The main contributions of our study comprise (i) a unified safeguarded framework for six classical \texttt{CG} methods; (ii) dimension-independent online coefficient learning for recovery and stagnation; (iii) independent online prediction for bounded line-search fallback corrections; (iv) improved robustness for all six variants, with broadly comparable evaluation costs and higher runtimes.

The remainder of the study is organized as follows. \Cref{sec:preliminaries} introduces the classical \texttt{HS}, \texttt{PRP}, \texttt{HZ}, \texttt{DY}, \texttt{FR}, and \texttt{DL} methods, together with a reformulation of the basic descent conditions used throughout the work. \Cref{sec:stagnation} describes the persistent-stagnation criterion and the dimension-independent feature representation. \Cref{sec:nn-parameters} presents the online \texttt{ELM} coefficient-prediction mechanism, the recovery procedure, the stagnation intervention, and the associated deterministic safeguards. \Cref{sec:training} describes the online training procedure, the data archive, and the model-update strategy, and introduces the learning-assisted fallback-step mechanism within the existing \texttt{cvsrch} failure branch. \Cref{sec:computational-cost} 
analyzes the additional computational cost imposed by the learning-assisted mechanism relative to the corresponding classical \texttt{CG} methods, together with discussing the global convergence. \Cref{sec:experiments} presents the computational experiments and compares the proposed classical methods with their corresponding learning-assisted variants. Finally, \Cref{sec:conclusion} provides the concluding remarks.

\section{Preliminaries}
	\label{sec:preliminaries}

This section presents the main analytical components underlying the proposed learning-assisted \texttt{CG} framework. We first describe the line-search strategy used to determine the step length and recall the \texttt{SDC} imposed on the search directions. We then review the traditional \texttt{CG} coefficients employed as the analytical basis of the proposed algorithms, together with the corresponding restart strategies used when these coefficients become unreliable or unavailable. These classical ingredients provide the reference mechanism against which the neural interventions are subsequently defined and safeguarded. In particular, the learning component is designed to complement, rather than systematically replace, the underlying analytical \texttt{CG} scheme.

\subsection{Line-Search Strategy}
	
As is well known, when a descent direction is given, a line-search procedure determines an appropriate step length along that direction. On iterations for which the \texttt{cvsrch} procedure succeeds, the step length in \eqref{eq:iteration} is required to satisfy the well-known Wolfe conditions \cite{nocedal1999numerical}; that is,
\begin{equation}
\label{eq:wolfe-armijo}
f(\x_k+\alpha_k\dvec_k) \leq f(\x_k)+c_1\alpha_k\g_k^{\mathsf T}\dvec_k,
\end{equation}
and
\begin{equation}
\label{eq:wolfe-curvature}
\nabla f(\x_k+\alpha_k\dvec_k)^{\mathsf T}\dvec_k \geq c_2\g_k^{\mathsf T}\dvec_k,
\end{equation}
implemented using the {\tt cvsrch} solver \cite{cvsrch}, where \(0<c_1<c_2<1\) are prespecified constants. As a stronger requirement than the standard descent property, the search direction is often required to satisfy the \texttt{SDC}. Let
\begin{equation}
\label{eq:ck-Gk}
c_k:=\g_k^{\mathsf T}\dvec_{k-1},
\qquad
G_k:=\norm{\g_k}^{2},
\qquad \text{and}\qquad
D_{k-1}:=\norm{\dvec_{k-1}}^{2}.
\end{equation}
Considering the two-term \texttt{CG} form of the search direction in
\eqref{eq:two-term-direction}, the \texttt{SDC} is equivalent to
\begin{equation}
\label{eq:sdc-beta}
\beta_k c_k\leq(1-\delta)G_k.
\end{equation}
Thus, a coefficient \(\beta_k\) is compatible with the \texttt{SDC}
whenever \eqref{eq:sdc-beta} holds.

\paragraph{{\bf A uniform safeguard for the search direction}}
We impose the \texttt{SDC} and \texttt{DSC} defined in
\eqref{eq:ls-m3-direction-bounds}. Within each learning-assisted variant,
both conditions are verified before calling \texttt{cvsrch}. By the two-term formula
\eqref{eq:two-term-direction}, the \texttt{DSC} is equivalent to
\begin{equation}
\label{eq:ls-m3-scale-quadratic}
D_{k-1}\beta_k^2-2c_k\beta_k+(1-\kappa_d^2)G_k\leq0.
\end{equation}
Consequently, both bounds hold whenever
\begin{equation}
\label{eq:ls-m3-admissible}
\beta_k\in\mathcal K_k
:=\bigl\{\beta\in\mathbb R:\;
\beta c_k\leq(1-\delta)G_k,\ \ 
D_{k-1}\beta^2-2c_k\beta+(1-\kappa_d^2)G_k\leq0\bigr\}.
\end{equation}
If \(D_{k-1}>0\), the set \(\mathcal K_k\) is a nonempty compact
interval because \(0\in\mathcal K_k\) for \(\kappa_d\geq1\).
If \(D_{k-1}=0\), then \(\dvec_{k-1}=0\), \(c_k=0\),
and the direction equals \(-\g_k\) for every coefficient; in this
case we simply choose \(\beta_k=0\). For \(D_{k-1}>0\), the scale constraint in
\eqref{eq:ls-m3-admissible} is the interval
\begin{equation}
	\label{eq:ls-m3-scale-interval}
	\left[
	\frac{c_k-\sqrt{c_k^2+(\kappa_d^2-1)D_{k-1}G_k}}{D_{k-1}},
	\frac{c_k+\sqrt{c_k^2+(\kappa_d^2-1)D_{k-1}G_k}}{D_{k-1}}
	\right].
\end{equation}
For \(D_{k-1}>0\), given a finite proposed coefficient
\(\beta_k^{\mathrm{prop}}\), it is retained exactly if it belongs to
\(\mathcal K_k\); otherwise the Euclidean projection
\begin{equation}
\label{eq:ls-m3-projection}
\beta_k^{\mathrm{used}}
=\Pi_{\mathcal K_k}\bigl(\beta_k^{\mathrm{prop}}\bigr),\qquad
\dvec_k=-\g_k+\beta_k^{\mathrm{used}}\dvec_{k-1},
\end{equation}
is selected. At the final direction-guard stage, a nonfinite selected coefficient,
an unusable previous direction, or failure of numerical verification
leads to the \texttt{SD} direction defined in \eqref{eq:steepest-descent-direction}; when \(\g_k=0\), the method terminates. Method-specific restart rules
remain in force, and any subsequent adjustment of the coefficient is
followed by the same two tests. This retains the classical direction
unchanged whenever it already satisfies both requirements.

The following result formalizes the simultaneous enforcement of the two
conditions. It extends the sufficient-descent repair principle in
\cite[Proposition~3.1]{Kimiaei2026}, which for the Euclidean norm
repairs a tentative vector \(v\) by setting
\[
\widehat v=v-\left(\delta+{\g_k^{\mathsf T}v}/{G_k}\right)\g_k
\qquad \text{and}\qquad \g_k^{\mathsf T}\widehat v=-\delta G_k.
\]
This correction alone does not imply a bound on
\(\norm{\widehat v}/\norm{\g_k}\); our coefficient restriction
enforces both inequalities while retaining the two-term \texttt{CG} form.

\begin{theorem}[Simultaneous descent and gradient-relative direction bound]
\label{thm:ls-m3-simultaneous}
Let \(\g_k\in\R^n\setminus\{0\}\) and the preceding direction
\(\dvec_{k-1}\in\R^n\) be finite, and fix
\(0<\delta<1\) and \(\kappa_d\geq1\). Define \(c_k,G_k,D_{k-1}\)
and \(\mathcal K_k\) by \eqref{eq:ck-Gk} and
\eqref{eq:ls-m3-admissible}, and let \(\beta_k^{\mathrm{prop}}\in\R\)
be any finite proposed coefficient. If
\(D_{k-1}>0\), choose \(\beta_k^{\mathrm{used}}\) by
\eqref{eq:ls-m3-projection}; otherwise set \(\beta_k^{\mathrm{used}}=0\).
Then the resulting two-term direction
\(\dvec_k=-\g_k+\beta_k^{\mathrm{used}}\dvec_{k-1}\) satisfies the
\texttt{SDC} and \texttt{DSC} in
\eqref{eq:ls-m3-direction-bounds}, and additionally
\begin{equation}
\label{eq:lowdirnorm}
\norm{\dvec_k}\ge \delta\norm{\g_k}.
\end{equation}
In particular, its descent-angle cosine satisfies
\begin{equation}
\label{eq:ls-m3-angle}
\frac{-\g_k^{\mathsf T}\dvec_k}
{\norm{\g_k}\norm{\dvec_k}}
\geq\frac{\delta}{\kappa_d}>0.
\end{equation}
If \(D_{k-1}>0\) and the proposed coefficient belongs to
\(\mathcal K_k\), it is retained exactly.
\end{theorem}

\begin{proof}
For \(D_{k-1}>0\), the scale constraint defining \(\mathcal K_k\)
is a convex quadratic inequality with positive leading coefficient;
its solution set is the compact interval displayed above. Moreover,
\(\beta=0\) satisfies this inequality because
\((1-\kappa_d^2)G_k\leq0\), and it satisfies the descent inequality
because \(0\leq(1-\delta)G_k\). Hence \(\mathcal K_k\) is
nonempty, closed, and convex. The Euclidean projection in
\eqref{eq:ls-m3-projection} is therefore well-defined and lies in
\(\mathcal K_k\).
Using the two-term formula gives
\[
\g_k^{\mathsf T}\dvec_k=-G_k+\beta_k^{\mathrm{used}}c_k\quad\text{and}\quad
\norm{\dvec_k}^2=G_k-2\beta_k^{\mathrm{used}}c_k
+(\beta_k^{\mathrm{used}})^2D_{k-1}.
\]
The two defining inequalities of \(\mathcal K_k\) now yield the \texttt{SDC} and \texttt{DSC} in \eqref{eq:ls-m3-direction-bounds}, i.e., \(\g_k^{\mathsf T}\dvec_k\leq-\delta G_k\) and
\(\norm{\dvec_k}^2\leq\kappa_d^2G_k\). By the Cauchy--Schwarz inequality,
\(\delta G_k\leq-\g_k^{\mathsf T}\dvec_k
\leq\norm{\g_k}\norm{\dvec_k}\), which yields the lower norm
bound \eqref{eq:lowdirnorm}. Combining the descent and upper norm inequalities proves \eqref{eq:ls-m3-angle}. If \(D_{k-1}=0\), the choice \(\beta_k=0\) yields the \texttt{SD} direction in \eqref{eq:steepest-descent-direction}, which satisfies every conclusion. Finally, the projection of a point already in the admissible set is the point itself.
\end{proof}

The theorem applies to the final selected direction at every
nonstationary iteration in the exact-arithmetic model of each
learning-assisted variant. The implementation separately verifies
finiteness and its prescribed numerical descent and direction-scale
tests, using the \texttt{SD} direction defined in
\eqref{eq:steepest-descent-direction} whenever final verification fails.

\subsection{Traditional Conjugate Gradient Coefficients}

Let the step \(\svec_{k}\) and the corresponding gradient displacement \(\yvec_k\) be defined by
\begin{equation}
	\label{eq:secant-vectors}
	\svec_{k}:=\x_k-\x_{k-1} \quad\text{and}\quad \yvec_k:=\g_k-\g_{k-1}.
\end{equation}
The learning-assisted framework considered in this work is built upon the following basic \texttt{CG} coefficient formulas:
\begin{equation}
	\label{eq:cg-parameters}
		\beta_k^{\mathrm{PRP}}
		=
		\frac{\g_k^{\mathsf T}\yvec_k}
		{\norm{\g_{k-1}}^2}, \ \ \beta_k^{\mathrm{HS}}
		=
		\frac{\g_k^{\mathsf T}\yvec_k}
		{\dvec_{k-1}^{\mathsf T}\yvec_k}, \  \
		\beta_k^{\mathrm{FR}}
		=
		\frac{\norm{\g_k}^2}
		{\norm{\g_{k-1}}^2}, \  \
        \beta_k^{\mathrm{DY}}
		=
		\frac{\norm{\g_k}^2}
		{\dvec_{k-1}^{\mathsf T}\yvec_k}.
\end{equation}
Moreover, our study also includes the well-known \texttt{DL} method \cite{SBKRairo}, which incorporates an additional scalar \(t_k>0\) into the \texttt{CG} coefficient as follows: 
\begin{equation}
	\label{eq:dl-beta-introduction}
	\beta_k^{\mathrm{DL}}(t_k)
	:=
	\frac{\g_k^{\mathsf T}\yvec_k}{\dvec_{k-1}^{\mathsf T}\yvec_k}
	-
	t_k\frac{\g_k^{\mathsf T}\svec_{k}}{\dvec_{k-1}^{\mathsf T}\yvec_k}.
\end{equation}
The parameter \(t_k\) is not predicted by the learning model. For
\texttt{DL} and \texttt{DLNN}, we use the fixed value $t_k=t_{\mathrm{DL}}>0$ (cf. \cite{DaiLiaoNCG}). By contrast, \texttt{HZ} and \texttt{HZNN} use
\begin{equation}
    \label{eq:dl-t-introduction}
    t_k^{\mathrm{HZ}}
    :=
    2\frac{\norm{\yvec_k}^{2}}
    {\svec_{k}^{\mathsf T}\yvec_k},
\end{equation}
whenever the required quantities are well-defined (cf. \cite{hager2005new}). Thus, the \texttt{HZ} coefficient is the \texttt{DL} expression \eqref{eq:dl-beta-introduction} evaluated with
\(t_k=t_k^{\mathrm{HZ}}\), whereas the \texttt{DL} coefficient uses the
fixed parameter \(t_{\mathrm{DL}}\). The same respective choices are used
in \texttt{HZNN} and \texttt{DLNN}.

At an ordinary iteration, the selected analytical coefficient
\(\beta_k^{\mathrm{raw}}\) is substituted into the two-term search direction update
\eqref{eq:two-term-direction}. Hence,
\begin{equation}
	\label{eq:classical-cg-direction}
	\dvec_k^{\mathrm{class}}
	:=
	-\g_k
	+
	\beta_k^{\mathrm{raw}}\dvec_{k-1},
\end{equation}
where
\begin{equation}
	\label{eq:raw-beta-method}
	\beta_k^{\mathrm{raw}}
	\in
	\left\{
	\beta_k^{\mathrm{PRP}},
    \beta_k^{\mathrm{HS}},
	\beta_k^{\mathrm{FR}},
	\beta_k^{\mathrm{DY}},
	\beta_k^{\mathrm{DL}},
	\beta_k^{\mathrm{HZ}}
	\right\}.
\end{equation}
According to the solver state, the analytical coefficient remains the
primary coefficient in every learning-assisted variant. Within such a
variant, if \(\beta_k^{\mathrm{raw}}\) is finite and nonzero,
no method-specific restart condition is active, and the
corresponding analytical direction
\(\dvec_k^{\mathrm{class}}\) satisfies the \texttt{SDC} and
\texttt{DSC} in \eqref{eq:ls-m3-direction-bounds}, then it is retained
exactly, except when the separate persistent-stagnation mechanism becomes
eligible. If \(\beta_k^{\mathrm{raw}}\) and the corresponding analytical
direction are finite but the direction violates either the \texttt{SDC}
or the \texttt{DSC}, the analytical coefficient is deterministically
projected onto the admissible set \(\mathcal K_k\) in
\eqref{eq:ls-m3-admissible} before \texttt{cvsrch} is called. This projection is applied in the learning-assisted variants and does not
modify the corresponding standalone traditional \texttt{CG} methods.

To mitigate the negative effects of jamming \cite{hager2006survey}, the \texttt{DY} and \texttt{FR} implementations here employ the additional restart condition
\begin{equation}
	\label{eq:dy-fr-restart}
	\left|\g_k^{\mathsf T}\yvec_k\right|
	\leq
	\sqrt{\epsilon_{\mathrm{mach}}}\,\norm{\g_k}^{2},
\end{equation}
where \(\epsilon_{\mathrm{mach}}\) denotes machine precision. When \eqref{eq:dy-fr-restart} is satisfied, the classical
\texttt{DY} or \texttt{FR} method uses the \texttt{SD} direction defined in \eqref{eq:steepest-descent-direction} \cite{BabaieKafaki2024}. In contrast, the corresponding \texttt{DYNN} and \texttt{FRNN} variants first invoke a new safeguarded recovery mechanism, together with the necessary scale restriction, as is discussed in \Cref{sec:nn-parameters}. Meanwhile, for \texttt{HS} and \texttt{PRP}, a nonfinite analytical coefficient causes the corresponding classical implementation to restart. For \texttt{HZ}, the coefficient is regarded as unavailable whenever the quantities required in \eqref{eq:dl-beta-introduction} or \eqref{eq:dl-t-introduction} cannot be computed finitely. Similar strategies are generally applied to \texttt{DL} as well. In the corresponding \texttt{NN} variants, such issues are handled by the common restart-recovery hierarchy rather than by introducing an alternative analytical reference formula.
%developed in \Cref{sec:nn-parameters}; the additional scale restriction associated with this restart event is introduced there.
% \begin{equation}
% 	\label{eq:classical-descent}
% 	\g_k^{\mathsf T}\dvec_k^{\mathrm{class}}<0,
% \end{equation}

\subsection{Neural-Network Preliminaries}
\label{subsec:nn-preliminaries}

The learning components employed in the present work are scalar-output {\tt ELMs} \cite{SunELM2014,Wang2022ELM}. Their role is deliberately
restricted: They do not generate iterates, coordinate-wise scaling matrices,
or complete search directions. Instead, the first \texttt{ELM} predicts a
scalar \texttt{CG} coefficient when coefficient recovery or persistent
stagnation makes a neural proposal eligible, whereas a second, independent
\texttt{ELM} predicts a scalar correction to the fallback step length
following failure of the existing \texttt{cvsrch} line search
\cite{cvsrch}. Unlike a fully trained multilayer {\tt NN}, an
\texttt{ELM} employs a single hidden layer whose input weights and biases are
generated once and subsequently kept fixed. Only the output weights are
fitted using the data collected during the current optimization run
\cite{SunELM2014,Wang2022ELM}.

For \(v,a,b\in\R\) with \(a\leq b\), we define
\begin{equation}
	\label{eq:clip-definition}
	\operatorname{clip}(v,a,b)
	:=
	\min\{b,\max\{a,v\}\}.
\end{equation}

Let \(\mathcal N_{\ba}:\R^{p}\rightarrow\R\) denote a scalar-output
\texttt{ELM} with \(h\) hidden units. For an input feature vector
\(\bu\in\R^{p}\), let \(\bmu\in\R^p\) and \(\bsigma\in\R^p\) denote the
componentwise sample mean and standard deviation associated with the
current fitted model. The standardized feature vector is
\begin{equation}
	\label{eq:elm-standardization}
	\bz := (\bu-\bmu)\oslash\bsigma,
\end{equation}
where $\oslash$ denotes componentwise division. To avoid division by zero or excessive amplification of numerical errors, let \(\sigma_{\min}>0\) denote the standard-deviation safeguard threshold. Each component \(\sigma_i\) that is nonfinite or satisfies \(\sigma_i<\sigma_{\min}\) is replaced by \(1\)
before the standardization in \eqref{eq:elm-standardization} is performed. The standardized components are subsequently clipped to a bounded interval as follows: 
\begin{equation}
	\label{eq:elm-feature-clipping}
	[\bz]_i \leftarrow
	\operatorname{clip}([\bz]_i,-c_{\mathrm{feat}},c_{\mathrm{feat}}),
	\qquad i=1,\ldots,p,
\end{equation}
where \(c_{\mathrm{feat}}\) is a positive constant. Let
\(
\bW\in\R^{p\times h}\) and \(
\bb\in\R^h\) denote the fixed random input weights and hidden-layer biases.  The hidden
representation is defined here by
\begin{equation}
	\label{eq:elm-hidden-layer}
	\bh(\bu)
	:=
	\tanh
	\left(
	\bW^{\mathsf T}\bz+\bb
	\right),
\end{equation}
where the hyperbolic tangent is applied componentwise \cite{Wang2022ELM}.  After appending a
constant output-bias component, the scalar \texttt{ELM} output is then
\begin{equation}
	\label{eq:elm-output}
	\mathcal N_{\ba}(\bu)
	=
	\ba^{\mathsf T}\bar{\bh}(\bu), \ \ \text{with} \ \ 	\bar{\bh}(\bu)
	:=
	\begin{bmatrix}
		\bh(\bu)\\
		1
	\end{bmatrix}
	\in\R^{h+1},
\end{equation}
where
\(\ba\in\R^{h+1}\) contains the trainable output weights and output bias. Suppose that the current replay archive contains \(N\) input-target pairs
\begin{equation*}
	\label{eq:elm-training-set}
	\mathcal T_N
	=
	\left\{
	(\bu_j,r_j)
	\right\}_{j=1}^{N}.
\end{equation*}
After standardization and evaluation of the fixed hidden layer, let
\begin{equation}
	\label{eq:elm-hidden-matrix}
	\bH
	=
	\begin{bmatrix}
		\bar{\bh}(\bu_1)^{\mathsf T}\\
		\vdots\\
		\bar{\bh}(\bu_N)^{\mathsf T}
	\end{bmatrix}
	\quad\text{and}\quad
	\br
	=
	\begin{bmatrix}
		r_1\\
		\vdots\\
		r_N
	\end{bmatrix}.
\end{equation}
The output coefficients are obtained from ridge-regularized least squares
\cite{HoerlKennard1970}; that is,
\begin{equation}
	\label{eq:elm-ridge-training}
	\ba
	=
	\left(
	\bH^{\mathsf T}\bH
	+
	\lambda\bR
	\right)^{-1}
	\bH^{\mathsf T}\br,
\end{equation}
where
\begin{equation}
	\label{eq:elm-ridge-matrix}
	\bR
	=
	\operatorname{diag}(1,\ldots,1,0).
\end{equation}
Here, \(\lambda>0\) is the ridge-regularization parameter. The final zero in \(\bR\) ensures that the appended output bias remains
unregularized. Thus, the hidden-layer parameters remain fixed throughout the optimization run, while only \(\ba\) is updated as new training samples become available.

The \texttt{OOD} detection is commonly used to identify
inputs that depart from the distribution represented by the training data,
thereby providing a mechanism for rejecting predictions in regions where the
learned model may be less reliable \cite{Yang2024OOD}.  Existing
\texttt{OOD} approaches include classification-, density-, distance-, and
reconstruction-based methods, with distance-based techniques frequently
assessing departures in a learned or constructed feature space.  Motivated by
this general principle, the present framework employs a lightweight
standardized-feature \texttt{OOD} screening rule to determine whether an
online \texttt{ELM} prediction is sufficiently supported by its current
training archive.  The specific screening criterion used here is introduced
later.

In the proposed optimization framework, this generic \texttt{ELM}
construction is instantiated in two independent models with distinct roles.
The coefficient-\texttt{ELM} operates on a dimension-independent
representation of the current \texttt{CG} state and provides a candidate
coefficient only when restart recovery or persistent stagnation makes a neural
intervention eligible. It is trained online using valid classical
\texttt{CG} coefficients generated during the optimization run. The second
model, the fallback-step-\texttt{ELM}, is associated exclusively with failure
of the existing \texttt{cvsrch} procedure and provides a safeguarded
correction to the original deterministic fallback proposal; its training
information is obtained from accepted points returned by \texttt{cvsrch},
including certified cached Armijo points. The two
models maintain separate feature representations, replay archives,
standardization statistics, and network parameters, and neither requires an
offline data set or additional function or gradient evaluations for training.
In both cases, neural outputs are treated only as proposals.
A prediction is accepted only if the model is available, the current
standardized feature vector passes the \texttt{OOD} screening,
and the corresponding deterministic safeguards are satisfied. If a proposal is unavailable or rejected, the
appropriate classical direction or deterministic fallback mechanism is
retained. The detailed feature constructions, training targets, model
parameters, prediction rules, \texttt{OOD} screening criterion, and safeguards
for the two \texttt{ELM}s are presented in the subsequent sections.

\section{Stagnation Detection and Neural Inputs}
\label{sec:stagnation}

In this section, we briefly discuss the soft computing aspects of our study and address the intervention of the previously mentioned {\tt NN} model in our algorithmic objectives. To this end, we first introduce our stagnation detection strategy.

\subsection{Activation Criterion}

The learning-assisted method incorporates a restart-recovery mechanism and a separate stagnation mechanism. The restart-recovery mechanism may be invoked whenever the analytical \texttt{CG} update is unusable, independently of whether stagnation has occurred. By contrast, a stagnation-based neural proposal is considered only when a persistent lack of progress is detected in both the objective value and the gradient norm.

For the stagnation test, the algorithm stores the \(q\) most recent objective values and Euclidean gradient norms. Define
\[
\mathcal{F}_k:=\left\{f_{k-q+1},\ldots,f_k\right\}
\ \ 
\text{and} \ \ 
\mathcal{G}_k:=\left\{\norm{\g_{k-q+1}},\ldots,\norm{\g_k}\right\},
\]
whenever \(q\) observations are available. The corresponding relative variation measures are defined as
\begin{equation}
	\label{eq:relative-f-variation}
	R_k^f:={\left(\displaystyle\max_{f_j\in\mathcal{F}_k}f_j-\displaystyle\min_{f_j\in\mathcal{F}_k}f_j\right)}\big/{\max\left\{1,\displaystyle\max_{f_j\in\mathcal{F}_k}|f_j|\right\}}
\end{equation}
and
\begin{equation}
	\label{eq:relative-g-variation}
	R_k^g:={\left(\displaystyle\max_{\norm{\g_j}\in\mathcal{G}_k}\norm{\g_j}-\displaystyle\min_{\norm{\g_j}\in\mathcal{G}_k}\norm{\g_j}\right)}\big/{\max\left\{\displaystyle\max_{\norm{\g_j}\in\mathcal{G}_k}\norm{\g_j},\epsilon_{\mathrm{mach}}\right\}}
\end{equation}
with the machine precision \(\epsilon_{\mathrm{mach}}\). Thus, the detector measures the variation over the recent trajectory rather than the decrease between two consecutive iterates. This is particularly useful because neither the gradient norm nor its iteration-to-iteration variation is necessarily monotone.

Persistent stagnation requires, for prescribed positive constants
\(\tau_f\) and \(\tau_g\), both \(R_k^f\leq\tau_f\) and \(R_k^g\leq\tau_g\). A small variation in the
objective alone does not activate the stagnation mechanism, nor does a small
variation in the gradient norm alone. Furthermore, a stagnation-based neural
attempt is allowed only after a minimum number of iterations and only if
sufficiently many iterations have elapsed since the previous stagnation
attempt. Let \(k_{\min}\) denote the minimum activation iteration and
\(q_{\mathrm{gap}}\) the required spacing between consecutive stagnation
attempts. The complete stagnation trigger is therefore
\begin{equation}
	\label{eq:activation}
	\mathcal{A}_k:=\left\{R_k^f\leq\tau_f\right\}\cap\left\{R_k^g\leq\tau_g\right\}\cap\left\{k\geq k_{\min}\right\}\cap\left\{k-k_{\mathrm{last}}\geq q_{\mathrm{gap}}\right\},
\end{equation}
where \(k_{\mathrm{last}}\) denotes the iteration (index) of the most recent stagnation attempt. Condition \eqref{eq:activation} makes the neural model eligible for a stagnation-based intervention; it does not, by itself, replace the current direction. First, the ordinary analytical \texttt{CG} update and the restart-recovery mechanism are processed. The stagnation \texttt{ELM} is considered only if the resulting direction remains the valid classical direction. Accordingly, when persistent stagnation is detected and the safeguarded
\texttt{ELM} proposal is admissible, the method sets
\(\dvec_k=\dvec_k^{\texttt{NN}}\); otherwise, it retains
\(\dvec_k=\dvec_k^{\mathrm{class}}\). Here,
\(\dvec_k^{\mathrm{class}}\) denotes the valid classical \texttt{CG}
direction, while \(\dvec_k^{\texttt{NN}}\) denotes the corresponding
safeguarded \texttt{ELM}-based alternative. Hence, persistent stagnation makes the neural proposal eligible for consideration but does not imply its automatic acceptance. If the \texttt{ELM} is unavailable, the current feature vector is classified as \texttt{OOD}, or the resulting direction fails the prescribed safeguarding conditions, the classical direction is retained.
% In the implementation,
% \begin{equation*}
% 	\label{eq:stagnation-parameters}
% 	q=5, \qquad \tau_f=10^{-3}, \qquad \tau_g=10^{-3}.
% \end{equation*}
%The implemented values are \(k_{\min}=20\) and \(q_{\mathrm{gap}}=5\).

\subsection{Dimension-Independent Coefficient Feature Vector}

The coefficient-\texttt{ELM} predicts a single scalar \texttt{CG} coefficient. Accordingly, its input is constructed from scalar quantities rather than from the components of the \(n\)-dimensional optimization vectors. As a result, the input dimension is independent of the problem dimension.

For two vectors \(\ba\) and \(\bb\), define the regularized cosine quantity
\begin{equation}
	\label{eq:feature-cosine}
	C(\ba,\bb):=\frac{\ba^{\mathsf T}\bb}{\norm{\ba}\norm{\bb}+\epsilon_{\mathrm{mach}}}
\end{equation}
with the machine precision \(\epsilon_{\mathrm{mach}}\). Also, let the signed logarithmic transformation be defined as
\begin{equation}
	\label{eq:signed-log}
	L(v):=\operatorname{sign}(v)\log(1+|v|).
\end{equation}
Recalling \(\beta_k^{\mathrm{raw}}\) from \eqref{eq:raw-beta-method} and \(\dvec_k^{\mathrm{raw}}=\dvec_k^{\mathrm{class}}\) from \eqref{eq:classical-cg-direction}, $\chi_k^{\mathrm{raw}}:=\g_k^{\mathsf T}\dvec_k^{\mathrm{raw}}$. Then, the twelve-dimensional feature vector used by the coefficient-\texttt{ELM} is defined as follows:
\begin{equation}
	\label{eq:nn-input}
	\begin{aligned}
	\bu_k^{\beta}:=\Bigg[&
	\log\!\left(\frac{\norm{\g_k}+\epsilon_{\mathrm{mach}}}{\norm{\g_{k-1}}+\epsilon_{\mathrm{mach}}}\right),\,
	C(\g_k,\g_{k-1}),\,
	C(\g_k,\yvec_k),\,
	C(\g_k,\dvec_{k-1}),\\
	&
	C(\dvec_{k-1},\yvec_k),\,
	L(\dvec_{k-1}^{\mathsf T}\yvec_k),\,
	L(\g_k^{\mathsf T}\dvec_{k-1}),\,
	L(\svec_k^{\mathsf T}\yvec_k),\\
	&
	L(\beta_k^{\mathrm{raw}}),\,
	L(\chi_k^{\mathrm{raw}}),\,
	L(\beta_{k-1}^{\mathrm{used}}),\,
	\gamma_{\texttt{CG}}
	\Bigg]^{\mathsf T}.
	\end{aligned}
\end{equation}
Here, \(\beta_{k-1}^{\mathrm{used}}\) denotes the most recently used coefficient stored by the learning-assisted solver, while \(\gamma_{\texttt{CG}}\) identifies the corresponding analytical \texttt{CG} coefficient formula:
\begin{equation*}
	\label{eq:method-code}
	\gamma_{\texttt{CG}} \in
\left\{
\frac{2}{6}\;(\texttt{PRP}),\;
\frac{1}{6}\;(\texttt{HS}),\;
\frac{5}{6}\;(\texttt{FR}),\;
\frac{4}{6}\;(\texttt{DY}),\;
\frac{6}{6}\;(\texttt{DL}),\;
\frac{3}{6}\;(\texttt{HZ})
\right\}.
\end{equation*}
The feature vector in \eqref{eq:nn-input} contains no raw \(n\)-dimensional vectors and therefore has the fixed dimension \(\dim\left(\bu_k^{\beta}\right)=12\), independently of the model dimension. Before being supplied to the \texttt{ELM}, any nonfinite component of the
feature vector is replaced by zero, and each component is clipped using
the rule in \eqref{eq:elm-feature-clipping}, with \(\bz\) replaced by
\(\bu_k^{\beta}\) and \(p=12\). During model fitting, the archived feature
vectors are standardized using their sample means and standard deviations
before being passed through the fixed hidden layer of the \texttt{ELM}.

It is important to distinguish the coefficient history from the neural input. The implementation stores a bounded recent history of accepted finite coefficients; however, this history is not concatenated with \(\bu_k^{\beta}\), and neither forgetting factors nor mask variables are used. The stored coefficient history serves only as a deterministic bootstrap source when an \texttt{ELM}-based recovery coefficient is unavailable. Likewise, the stagnation measures \(R_k^f\) and \(R_k^g\) are not components of the coefficient feature vector in \eqref{eq:nn-input}; instead, they serve as the activation gate for the stagnation mechanism. These measures are included in the input of the separate line-search fallback \texttt{ELM} discussed later, but not in the twelve-dimensional feature vector used by the coefficient-\texttt{ELM}.

\section{Safeguarded Neural Parameter Prediction}
\label{sec:nn-parameters}

In this section, we describe how the \texttt{ELM} predictions are incorporated into the proposed learning-assisted optimization framework. Neural outputs are treated only as candidate parameters and are accepted only after the corresponding screening and deterministic safeguards are satisfied. We first address the safeguarded prediction of the \texttt{CG} coefficient and its use in recovery and stagnation mechanisms, and then discuss the independent step prediction employed after line-search failure.

\subsection{Online \texttt{ELM} Prediction of a \texttt{CG} Coefficient}

The coefficient-\texttt{ELM} uses the twelve-dimensional feature vector \(\bu_k^{\beta}\) defined by \eqref{eq:nn-input}. Feature standardization and clipping are performed according to \eqref{eq:elm-standardization}--\eqref{eq:elm-feature-clipping}, and the scalar \texttt{ELM} output is evaluated using \eqref{eq:elm-hidden-layer}--\eqref{eq:elm-output}. Moreover, the \texttt{ELM} output is mapped to the bounded coefficient as follows:
\begin{equation}
\label{eq:beta-elm-prediction}
\widetilde{\beta}_k=\beta_{\max}\tanh\!\left(\widehat r_k^\beta\right),
\end{equation}
where \(\beta_{\max}>0\) is a prescribed constant and
\(\widehat r_k^{\beta}=\mathcal N_{\ba_{\beta}}\left(\bu_k^{\beta}\right)\)
denotes the scalar output of the coefficient-\texttt{ELM}. Hence, the neural output is automatically bounded by
\(-\beta_{\max}<\widetilde{\beta}_k<\beta_{\max}\).
Equation \eqref{eq:beta-elm-prediction} represents an absolute coefficient prediction rather than a residual correction to the current analytical coefficient, which is adopted as the primary choice. If it is finite and nonzero, no method-specific restart condition
is active, and its corresponding classical \texttt{CG} direction satisfies the \texttt{SDC} and \texttt{DSC} in \eqref{eq:ls-m3-direction-bounds}, that direction is retained unless the persistent-stagnation mechanism described in \Cref{sec:stagnation} becomes active. Thus, the \texttt{ELM} prediction is used only as a candidate within specific recovery and stagnation procedures, as described below.

\subsection{{\tt OOD} Screening}

The {\tt OOD} test serves as a screening mechanism for the
coefficient-\texttt{ELM} introduced in the preceding subsection. Its purpose
is to prevent the neural coefficient \(\widetilde{\beta}_k\) from being used
when the current optimization state is too far from the states represented
in the coefficient-training archive. Thus, the {\tt OOD} test neither
generates a coefficient nor replaces the descent safeguard; it only
determines whether the coefficient-\texttt{ELM} is eligible to provide a
candidate coefficient. Motivated by the general principle of identifying
predictions whose feature representations depart from the distribution
represented by the training data \cite{Yang2024OOD}, we employ a simple
standardized-feature screening rule. The particular infinity-norm and
percentile-based criterion used here is introduced specifically as a
safeguard for the online coefficient-\texttt{ELM}.

To formulate this criterion, recall that the coefficient-\texttt{ELM}
receives the twelve-dimensional feature vector \(\bu_k^{\beta}\).
Let \(\boldsymbol{\mu}_{\beta}\) and \(\boldsymbol{\sigma}_{\beta}\)
denote the mean and standard deviation stored from the most recent
coefficient-\texttt{ELM} fit. Applying the generic standardization rule
already defined in \eqref{eq:elm-standardization} with these stored
statistics, let \(\zvec_j\) denote the standardized feature vector
associated with the \(j\)th archived training sample \(\bu_j^{\beta}\),
for \(j=1,\ldots,N_{\beta}^{\mathrm{fit}}\), where
\(N_{\beta}^{\mathrm{fit}}\) is the number of training samples used in
that fit. Likewise, let \(\zvec_k\) denote the standardized version of
the current feature vector \(\bu_k^{\beta}\), computed using the same
\(\boldsymbol{\mu}_{\beta}\) and \(\boldsymbol{\sigma}_{\beta}\).
Thus, \(j\) indexes the archived samples used to fit the current model,
whereas \(k\) denotes the current optimization iteration at which a
coefficient-\texttt{ELM} prediction is being considered.

The {\tt OOD} screening is based on the infinity norm of these standardized
feature vectors. During each coefficient-\texttt{ELM} fit, the standardized
archived feature vectors are clipped componentwise using the same feature
bound \(c_{\mathrm{feat}}>0\) employed by the model. Accordingly, for the
\(j\)th archived standardized feature vector \(\zvec_j\), define
\(s_j:=\min\{c_{\mathrm{feat}},\norm{\zvec_j}_{\infty}\}\), so that
\(s_j\) represents the largest absolute standardized component after the
feature clipping used in the implementation. Let
\(\nu_{\texttt{OOD}}\in(0,1)\) denote the archive quantile used in the
\texttt{OOD} scale. If
\(s_{(1)}\leq s_{(2)}\leq\cdots\leq
s_{(N_{\beta}^{\mathrm{fit}})}\) denote the sorted archive scores, define
\begin{equation*}
	Q_{\nu_{\texttt{OOD}}}\left(
	\{s_j\}_{j=1}^{N_{\beta}^{\mathrm{fit}}}
	\right)
	:=
	s_{\left(
		\left\lceil
		\nu_{\texttt{OOD}}N_{\beta}^{\mathrm{fit}}
		\right\rceil
		\right)}.
\end{equation*}
After each coefficient-\texttt{ELM} fit, the implementation therefore
summarizes the scale of the current training archive by
\begin{equation*}
	s_{{\tt OOD}}
	:=
	\max\left\{
	1,\,
	Q_{\nu_{\texttt{OOD}}}\left(
	\{s_j\}_{j=1}^{N_{\beta}^{\mathrm{fit}}}
	\right)
	\right\}.
\end{equation*}
Hence, \(s_{{\tt OOD}}\) provides a reference scale determined by the
standardized and clipped feature vectors represented in the archive used
to fit the current model.

Using the current standardized feature vector \(\zvec_k\) defined above,
and before applying the feature clipping used for \texttt{ELM} evaluation,
the current state is regarded as sufficiently close to the training
distribution only if
\begin{equation}
	\label{eq:ood-test}
	\norm{\zvec_k}_{\infty}
	\leq
	\gamma_{{\tt OOD}}\,s_{{\tt OOD}},
\end{equation}
where \(\gamma_{{\tt OOD}}>0\) is a prescribed constant. Thus,
\(\gamma_{{\tt OOD}}\) controls how far a current standardized feature
vector may extend beyond the scale represented by the clipped training
archive before the neural prediction is rejected. The role of
\eqref{eq:ood-test} in the coefficient-selection procedure is therefore
as follows. When the coefficient-\texttt{ELM} is invoked, it is not allowed
to provide the candidate coefficient \(\widetilde{\beta}_k\) if the model
has not yet been trained, if \(\zvec_k\) contains any nonfinite component,
or if \eqref{eq:ood-test} is not satisfied. Once these eligibility checks
are passed, \(\zvec_k\) is clipped componentwise to the interval
\([-c_{\mathrm{feat}},c_{\mathrm{feat}}]\) before being supplied to the
trained \texttt{ELM}, which generates the candidate coefficient
\(\widetilde{\beta}_k\) defined in \eqref{eq:beta-elm-prediction}.
Passing the {\tt OOD} test does not, by itself, imply that this coefficient
will be used; the resulting candidate must still satisfy the deterministic
admissibility and descent safeguards described in the next subsections.

The consequence of {\tt OOD} rejection depends on the reason for invoking
the coefficient-\texttt{ELM}. In the restart-recovery mechanism, the
classical coefficient or direction is already unusable; therefore,
rejection of the neural proposal causes the algorithm to proceed to the
coefficient-history and deterministic bootstrap hierarchy. In the
persistent-stagnation mechanism, by contrast, the classical direction has
already passed the uniform safeguard \eqref{eq:ls-m3-direction-bounds}.
Consequently, if the current feature vector fails the \texttt{OOD} test,
the neural alternative is simply rejected and the already-admissible
classical \texttt{CG} direction is retained. Hence, \(s_j\),
\(s_{{\tt OOD}}\), and \(\gamma_{{\tt OOD}}\) are not optimization parameters
that directly modify the \texttt{CG} formula. They are used solely to
determine whether the current coefficient-\texttt{ELM} prediction is
sufficiently supported by the feature distribution represented in the
corresponding online training archive.

\subsection{Sufficient Descent Safeguard}

The \texttt{SDC} and \texttt{DSC} have already been defined in \eqref{eq:ls-m3-direction-bounds}. For each of the six learning-assisted variants, the final search direction is explicitly checked before \texttt{cvsrch} is called. Thus, within the learning-assisted methods, the \texttt{SDC} and \texttt{DSC} are deterministic algorithmic safeguards
rather than assumptions imposed on the generated directions. The implementation applies these safeguards in two stages: local coefficient recovery enforces the corresponding local \texttt{SDC} together with a coefficient bound, whereas the common final guard enforces both the \texttt{SDC} and the finite-precision form of the \texttt{DSC}.

The implementation distinguishes the parameter used by the local
coefficient-recovery routines from that used by the final common direction
guard. Define
\begin{equation}
\label{eq:implemented-delta-local}
\delta_{\mathrm{loc}}
:=
\max\left\{
\left|\Delta_{\texttt{SDC}}\right|,
\epsilon_{\mathrm{mach}}
\right\}.
\end{equation}
The local coefficient intervals and descent-repair routines use
\(\delta_{\mathrm{loc}}\).

Before \texttt{cvsrch} is called, however, every direction generated by a
learning-assisted variant is subjected to a common final verification using
\begin{equation}
\label{eq:implemented-delta}
\delta
:=
\min\left\{
\delta_{\mathrm{loc}},
1-\epsilon_{\mathrm{mach}}
\right\}.
\end{equation}
Hence \(0<\delta<1\), as required by the theoretical safeguard.

To avoid modifying the local descent interval based on a numerically
negligible value of \(c_k\) defined in \eqref{eq:ck-Gk}, the implementation
introduces
\begin{equation*}
\label{eq:ck-tolerance}
\tau_k^c
:=
\epsilon_{\mathrm{mach}}
\max\left\{
1,\,
\norm{\g_k}\norm{\dvec_{k-1}}
\right\}.
\end{equation*}
Together with the coefficient bound
\(\beta\in[-\beta_{\max},\beta_{\max}]\), the tolerance-relaxed local screening interval is
\begin{equation}
\label{eq:local-beta-admissible-set}
\widetilde{\mathcal I}_k^\beta
:=
\left\{
\beta\in[-\beta_{\max},\beta_{\max}]
:
|c_k|\leq\tau_k^c
\quad\text{or}\quad
\beta c_k\leq(1-\delta_{\mathrm{loc}})G_k
\right\},
\end{equation}
implemented as
\begin{equation*}
\widetilde{\mathcal I}_k^\beta
=
\begin{cases}
\left[
-\beta_{\max},
\,
\min\left\{
\beta_{\max},
(1-\delta_{\mathrm{loc}})G_k/c_k
\right\}
\right],
& c_k>\tau_k^c,
\\[1ex]
\left[
\max\left\{
-\beta_{\max},
(1-\delta_{\mathrm{loc}})G_k/c_k
\right\},
\,
\beta_{\max}
\right],
& c_k<-\tau_k^c,
\\[1ex]
[-\beta_{\max},\,\beta_{\max}],
& |c_k|\leq\tau_k^c.
\end{cases}
\end{equation*}

For comparison with the final direction safeguard, define the exact-arithmetic
sufficient-descent set
\begin{equation}
\label{eq:beta-admissible-set}
\mathcal I_k^\beta
:=
\left\{
\beta\in[-\beta_{\max},\beta_{\max}]
:
\beta c_k
\leq
(1-\delta)G_k
\right\}.
\end{equation}
The direction-scale condition in
\eqref{eq:ls-m3-direction-bounds} gives
\begin{equation}
\label{eq:beta-scale-admissible-set}
\mathcal J_k^\beta
:=
\left\{
\beta\in\mathbb R:
D_{k-1}\beta^2
-2c_k\beta
+(1-\kappa_d^2)G_k
\leq0
\right\}.
\end{equation}
For \(D_{k-1}>0\), \(\mathcal J_k^\beta\) is the interval given in
\eqref{eq:ls-m3-scale-interval}. Since \(0\in\mathcal J_k^\beta\) for \(\kappa_d\geq1\), this interval is
nonempty. Recall $\mathcal K_k$ from \eqref{eq:ls-m3-admissible}; for recovery coefficients already restricted by
\([-\beta_{\max},\beta_{\max}]\), it is useful to define the complete
exact-arithmetic admissible set
\begin{equation}
\label{eq:complete-beta-admissible-set}
\widehat{\mathcal I}_k^\beta
:=
\mathcal I_k^\beta
\cap
\mathcal J_k^\beta
=
\mathcal K_k
\cap
[-\beta_{\max},\beta_{\max}].
\end{equation}
This set characterizes recovery coefficients satisfying both theoretical
direction safeguards. It should not, however, be interpreted as the interval
used by the initial local recovery stage.

Let \(\beta_k^{\mathrm{prop}}\) denote a coefficient supplied by the
coefficient-\texttt{ELM}, the recent coefficient history, or the
deterministic safe-bootstrap procedure. If \(\widetilde{\mathcal I}_k^\beta\) is empty, the recovery routine
returns the \texttt{SD} direction defined in
\eqref{eq:steepest-descent-direction}. Otherwise, the initial local
coefficient is obtained from
\begin{equation}
\label{eq:safe-beta}
\beta_k^{(0)}
:=
\Pi_{\widetilde{\mathcal I}_k^\beta}
\left(
\beta_k^{\mathrm{prop}}
\right).
\end{equation}
Thus, a finite proposal already belonging to
\(\widetilde{\mathcal I}_k^\beta\) is retained unchanged at this stage. If $|\beta_k^{(0)}|
\leq
\beta_{\min}^{\mathrm{nz}}$, the ordinary recovery routine attempts to construct a nonzero coefficient
inside \(\widetilde{\mathcal I}_k^\beta\), choosing its sign, when possible,
to favor descent. Hence \(\beta_{\min}^{\mathrm{nz}}\) is only a local
recovery threshold and is not a lower bound on the coefficient eventually
submitted to the line search.

If the resulting local direction fails the \texttt{SDC} in \eqref{eq:ls-m3-direction-bounds} with \(\delta=\delta_{\mathrm{loc}}\), the coefficient is successively contracted according to
$\beta_k^{(\ell+1)}
=
\vartheta_\beta\beta_k^{(\ell)}$, where $0<\vartheta_\beta<1$, and the direction is recomputed, up to the prescribed number of repair attempts. If no finite nonzero locally admissible coefficient can be constructed, the recovery routine returns the \texttt{SD} direction defined in \eqref{eq:steepest-descent-direction}.

After the analytical, recovery, or persistent-stagnation mechanism has
selected a tentative direction, the same common final guard is applied to
all learning-assisted variants. In particular, the final direction is
checked against the \texttt{SDC} and \texttt{DSC} defined in
\eqref{eq:ls-m3-direction-bounds}. To account for finite-precision effects,
let \(\chi_{\mathrm{fp}}>0\) denote the finite-precision inflation factor.
The implementation verifies the \texttt{DSC} using the slightly enlarged
constant $\kappa_d^{\mathrm{fp}}
:=
\kappa_d
\left(
1+\chi_{\mathrm{fp}}\epsilon_{\mathrm{mach}}
\right)$,  while the \texttt{SDC} is checked using the prescribed value of
\(\delta\). Thus, a tentative direction satisfying the \texttt{SDC} and the \texttt{DSC} (i.e., $\norm{\dvec_k}
\leq
\kappa_d^{\mathrm{fp}}
\norm{\g_k}$) is retained unchanged. Otherwise, the current coefficient is projected onto
the direction-based interval determined jointly by the \texttt{SDC} and
\texttt{DSC}, and the two-term direction is reconstructed. The
implementation performs this projection in a normalized form to avoid
squaring a potentially very large previous direction.

The reconstructed direction is then checked again against the
\texttt{SDC} and the finite-precision form of the \texttt{DSC} described
above. If it is nonfinite or still violates either safeguard, the method
uses the \texttt{SD} direction defined in
\eqref{eq:steepest-descent-direction}. Thus, every direction ultimately
supplied to \texttt{cvsrch} by a learning-assisted variant satisfies the
\texttt{SDC} and the finite-precision form of the \texttt{DSC}.

The coefficient cap \(\beta_{\max}\) belongs to the local neural and
bootstrap recovery mechanisms. It is not imposed as an independent final
restriction on every analytical coefficient. In particular, an analytical
coefficient whose corresponding direction already satisfies the
\texttt{SDC} and the finite-precision form of the \texttt{DSC} is retained
unchanged by the final direction guard. No corresponding common final
direction guard is imposed on the standalone traditional
\texttt{CG} implementations.

\subsection{Restart-Recovery Hierarchy}

Recall \(\beta_k^{\mathrm{raw}}\) from \eqref{eq:raw-beta-method} and
\(\dvec_k^{\mathrm{class}}\) from \eqref{eq:classical-cg-direction}.
Every finite proposed coefficient is first subject to the uniform
direction safeguard developed in \Cref{sec:preliminaries}. Hence, a finite
classical coefficient that violates one of the conditions in
\eqref{eq:ls-m3-direction-bounds} is adjusted by the projection
\eqref{eq:ls-m3-projection} before \texttt{cvsrch} is called.

The restart-recovery procedure is invoked when the analytical
\texttt{CG} update cannot provide a usable nonzero coefficient, namely,
when \(\beta_k^{\mathrm{raw}}\) is zero or nonfinite, when the corresponding
direction is nonfinite, or when a method-specific restart condition is
satisfied. In these situations, the learning-assisted method first attempts
to construct an admissible nonzero coefficient using the recovery hierarchy
described below before reverting to the \texttt{SD} direction defined in \eqref{eq:steepest-descent-direction}. 

The coefficient source is selected in the following order. First, if a trained \texttt{ELM} is available and \eqref{eq:ood-test} is satisfied, the method uses the neural proposal \eqref{eq:beta-elm-prediction}. If an admissible \texttt{ELM} prediction is unavailable, the method examines the recent history of finite coefficients. If this history is nonempty, it uses
\[
\beta_k^{\mathrm{hist}}:=\operatorname{median}\left(\mathcal{B}_k\right),
\]
where \(\mathcal{B}_k\) denotes the bounded recent history of stored finite classical coefficients. Meantime, if neither an \texttt{ELM} prediction nor a usable history value is available, a deterministic nonzero coefficient is constructed within the
local admissible interval \(\widetilde{\mathcal I}_k^\beta\) defined in
\eqref{eq:local-beta-admissible-set}. This safe-bootstrap coefficient is chosen, when possible, with a sign intended to oppose that of \(c_k\) defined in \eqref{eq:ck-Gk}. The resulting tentative direction is subsequently subjected to the common final direction guard and verified against both conditions in \eqref{eq:ls-m3-direction-bounds}. Thus, an unavailable or {\tt OOD} \texttt{ELM} does not automatically trigger a restart. A true restart occurs only if the \texttt{ELM}, coefficient-history bootstrap, and deterministic safe-bootstrap mechanisms all fail to produce an acceptable nonzero direction.

\subsection{Additional Safeguard for \texttt{DYNN} and \texttt{FRNN} Restarts}

The \texttt{DY} and \texttt{FR} methods employ the additional restart
criterion \eqref{eq:dy-fr-restart}, introduced in the Preliminaries. When
this condition is satisfied, the corresponding classical method replaces
the two-term \texttt{CG} direction with the \texttt{SD} direction defined in \eqref{eq:steepest-descent-direction}. By
contrast, the \texttt{DYNN} and \texttt{FRNN} variants permit a safeguarded
recovery attempt, but impose an additional restriction on the component
inherited from the previous direction.

Let \(\rho>0\) denote the maximum permitted ratio of the inherited
component to the current gradient norm, and define $\beta_k^{\max,\mathrm{restart}}
:=
\rho{\norm{\g_k}}/{\norm{\dvec_{k-1}}}$. 
The restart-specific inherited-component restriction is
\begin{equation}
\label{eq:restart-component-bound}
\norm{\beta_k\dvec_{k-1}}
\leq
\rho\norm{\g_k}.
\end{equation}
Equivalently, $|\beta_k|\leq\beta_k^{\max,\mathrm{restart}}$. Accordingly, the restart-specific interval is
\begin{equation}
\label{eq:restart-scale-interval}
\mathcal J_{k,\mathrm{restart}}^\beta
:=
\left[
-\beta_k^{\max,\mathrm{restart}},
\,
\beta_k^{\max,\mathrm{restart}}
\right],
\end{equation}
and the local coefficient used in this branch is restricted to
\begin{equation}
\label{eq:restart-admissible-intersection}
\widetilde{\mathcal I}_{k,\mathrm{restart}}^\beta
:=
\widetilde{\mathcal I}_k^\beta
\cap
\mathcal J_{k,\mathrm{restart}}^\beta,
\end{equation}
where \(\widetilde{\mathcal I}_k^\beta\) is the local sufficient-descent
interval defined by \eqref{eq:local-beta-admissible-set}. The restriction \eqref{eq:restart-component-bound}, together with the two-term direction \eqref{eq:two-term-direction}, ensures the \texttt{DSC} in \eqref{eq:ls-m3-direction-bounds} with \(\kappa_d=1+\rho\).

In finite precision, the restart-specific routine accepts its local
candidate only if it satisfies the \texttt{SDC} and \texttt{DSC} in
\eqref{eq:ls-m3-direction-bounds}, with \(\delta\) replaced by
\(\delta_{\mathrm{loc}}\) defined in
\eqref{eq:implemented-delta-local}, and \(\kappa_d\) replaced by
\((1+\rho)(1+\chi_{\mathrm{fp}}\epsilon_{\mathrm{mach}})\), respectively.

Unlike the ordinary recovery mechanism, no fixed lower bound such as
\(\beta_{\min}^{\mathrm{nz}}\) is imposed in this branch. When
\(\norm{\dvec_{k-1}}/\norm{\g_k}\) is large, the restart scale restriction
may legitimately require a coefficient much smaller than the ordinary
nonzero recovery threshold. Imposing that threshold here could therefore
destroy the scale safeguard.

If the restart-specific interval cannot be formed, the candidate
coefficient is nonfinite or zero, or the locally reconstructed direction
fails either the restart-specific \texttt{SDC} or \texttt{DSC}, the
recovery attempt is rejected and the method uses the \texttt{SD} direction
defined in \eqref{eq:steepest-descent-direction}.

After this restart-specific check, the resulting tentative direction is
still subjected to the common final direction guard used by all
learning-assisted variants. Consequently, the direction ultimately supplied
to \texttt{cvsrch} satisfies the global \texttt{SDC} and the
finite-precision form of the global \texttt{DSC} associated with
\eqref{eq:ls-m3-direction-bounds}.

\subsection{Persistent-Stagnation Prediction}

The persistent-stagnation mechanism uses the same coefficient-\texttt{ELM}
but follows a different fallback policy. It is considered only when the ordinary coefficient selector has retained the classical update unchanged, and no method-specific restart recovery is active.

If the stagnation trigger \eqref{eq:activation} is satisfied, a neural coefficient is considered only
when the coefficient model has been trained and the current feature vector
passes the \texttt{OOD} test \eqref{eq:ood-test}. Neither the
coefficient-history bootstrap nor the deterministic safe bootstrap is used
in this branch.

The neural coefficient is first passed through the local descent-recovery
procedure, producing a tentative direction
\(\widetilde{\dvec}_k^{\texttt{NN}}\). This direction is considered as a
stagnation alternative only when it differs sufficiently from the already
available classical candidate:
\begin{equation}
\label{eq:stagnation-direction-change}
{
\norm{
\widetilde{\dvec}_k^{\texttt{NN}}
-
\dvec_k^{\mathrm{class}}
}
}\big/{
\max\left\{
\norm{\dvec_k^{\mathrm{class}}},
\epsilon_{\mathrm{mach}}
\right\}
}
\geq
\eta_{\mathrm{dir}},
\end{equation}
where \(\eta_{\mathrm{dir}}>0\) is a prescribed constant. If the coefficient-\texttt{ELM} is unavailable, the current feature vector
is classified as out of distribution, the local descent recovery fails, or
\eqref{eq:stagnation-direction-change} is not satisfied, the classical
candidate is retained. The selected tentative direction, whether classical or neural, is then
passed through the common final direction guard before \texttt{cvsrch} is
called. Consequently, the final common projection may modify the tentative
direction and may reduce its eventual relative change below
\(\eta_{\mathrm{dir}}\).

\subsection{Line-Search Failure and Safeguarded Step Prediction}

The coefficient \texttt{ELM} does not replace the line search. After a
search direction has been selected, the solver calls the existing
\texttt{cvsrch} routine \cite{cvsrch}. If the point returned by
\texttt{cvsrch} cannot be certified for direct reuse, the search direction
is normalized as
\(\bar{\dvec}_k:=\dvec_k/\norm{\dvec_k}\), and
\(\phi_k:=\g_k^{\mathsf T}\bar{\dvec}_k\) is computed. The original solver
then determines its minimal fallback step \(a_k^{\min}\). As in \cite{LMBOPT}, in formulaic form,
\begin{equation}
	\label{eq:original-amin}
	a_k^{\min}=
	\begin{cases}
		\delta_a\left|f_k/\phi_k\right|, &\qquad \exists i:\ x_{k,i}=0,\;\bar d_{k,i}\neq0,\\[1ex]
		\delta_a\max\left\{\left|f_k/\phi_k\right|,\displaystyle\min_{i:\bar d_{k,i}\neq0}\left|x_{k,i}/\bar d_{k,i}\right|\right\}, &\qquad \exists i:\ \bar d_{k,i}\neq0,\\[1ex]
		1, &\qquad \bar{\dvec}_k=\mathbf{0},
	\end{cases}
\end{equation}
where \(\delta_a>0\) is a small positive constant.

\paragraph{{\bf Reuse of a returned line-search trial}}
Before initiating fresh fallback backtracking, the solver examines the point
already returned by \texttt{cvsrch}. A returned point lying on the current
search ray is retained directly when the line-search return certifies the
Wolfe conditions. A returned point that is not certified as a Wolfe point
may nevertheless be accepted as a cached Armijo fallback if it gives strict
decrease, satisfies the Armijo condition, and either $\alpha\geq\alpha_{\min}$, or satisfies the secant-type condition
\[
\kappa_d^2
\norm{
\g(\x_k+\alpha\dvec_k)-\g_k
}
\geq
2\vartheta(1-c_1)\delta
\norm{\dvec_k}.
\]
The latter condition provides a lower-step certificate under Lipschitz
continuity and therefore permits an already evaluated Armijo point with
\(\alpha<\alpha_{\min}\) to be reused without another objective or gradient
evaluation. If the returned point instead violates Armijo at a positive
step \(\alpha_k^{\mathrm{bad}}\), this rejected step may be used to reduce
the initial seed of the subsequent fresh backtracking sequence.

\paragraph{{\bf Safeguarded steps after line-search failure}}
The fallback value \(a_k^{\min}\) in \eqref{eq:original-amin} is only an
initial proposal and is never accepted automatically. Let $\ell_k:=\norm{\dvec_k}>0$ and $\bar{\dvec}_k:={\dvec_k}/{\ell_k}$, where $0<\alpha_{\min}\leq\alpha_{\max}$ and $0<\vartheta<1$. If \(a_k^{\min}\) is finite and positive, set
\(a_k^{\mathrm{sel}}:=a_k^{\min}\); otherwise set
\(a_k^{\mathrm{sel}}:=\ell_k\alpha_{\min}\). The repaired normalized
base step is
\begin{equation}
\label{eq:ls-m2-repaired-base}
\alpha_k^{\mathrm{base}}
:=
\operatorname{clip}\left(
{a_k^{\mathrm{sel}}}/{\ell_k},
\alpha_{\min},
\alpha_{\max}
\right),
\end{equation}
where $\operatorname{clip}$ is defined by \eqref{eq:clip-definition}.

For a learning-assisted state, let \(\widehat r_k^\alpha\) denote the
output of the fallback-step-\texttt{ELM}, and define
\begin{equation}
\label{eq:step-elm-correction}
C_k^\alpha
:=
\exp\!\left(\widehat r_k^\alpha\right).
\end{equation}
The multiplicative correction is restricted by
\begin{equation}
\label{eq:amin-correction-bounds}
1\leq C_{\min}\leq C_k^\alpha\leq C_{\max}.
\end{equation}
If the model is unavailable, the current feature vector fails the
\texttt{OOD} test \eqref{eq:ood-test}, the prediction is nonfinite, or
the correction is inadmissible, \(C_k^\alpha=1\) is used.

If the independent fallback-step-\texttt{ELM} supplies an admissible
multiplicative correction \(C_k^\alpha\), it is restricted by \eqref{eq:amin-correction-bounds}; otherwise \(C_k^\alpha=1\). The initial fallback trial is therefore
\begin{equation}
\label{eq:ls-m2-initial-step}
\alpha_k^{(0)}
:=
\operatorname{clip}\left(
C_k^\alpha\alpha_k^{\mathrm{base}},
\alpha_{\min},
\alpha_{\max}
\right).
\end{equation}

If the unsuccessful \texttt{cvsrch} call has already returned a positive
trial \(\alpha_k^{\mathrm{bad}}\) known to violate Armijo, the fresh
backtracking sequence may instead begin from
\begin{equation}
	\label{eq:ls-m2-reduced-seed}
	\min\left\{
	\alpha_k^{(0)},
	\vartheta\alpha_k^{\mathrm{bad}}
	\right\}.
\end{equation}
This avoids restarting above a smaller step that is already known to fail. A fallback trial is accepted only when it has a finite objective value,
gives strict decrease, and satisfies \eqref{eq:wolfe-armijo}. Otherwise, $\alpha\leftarrow\vartheta\alpha$ and the test is repeated. No Wolfe curvature condition is imposed on this fresh fallback-backtracking branch.

When separate objective and gradient routines are available, rejected
fallback trials require only objective evaluations, and the gradient is
evaluated after Armijo acceptance. When only a joint function-gradient
routine is available, each trial evaluates both quantities. Training and
prediction themselves require no additional objective or gradient
evaluations, whereas verification of a fallback trial may do so.

For completeness, if the gradient is \(L\)-Lipschitz along the fallback
trial segments, the descent lemma shows that
\eqref{eq:wolfe-armijo} holds whenever
\[
0<\alpha\leq \frac{2(1-c_1)(-\g_k^{\mathsf T}\dvec_k)}
{L\norm{\dvec_k}^2}
\quad\text{and hence whenever}\quad
0<\alpha\leq\frac{2(1-c_1)\delta}{L\kappa_d^2}.
\]
It follows from Theorem~\ref{thm:ls-m3-simultaneous} and the backtracking
rule that every mathematically accepted fallback step satisfies
\begin{equation}
\label{eq:ls-m2-uniform-decrease}
\alpha_k\geq\underline\alpha_{\mathrm F}
:=\min\!\left\{\alpha_{\min},
\frac{2\vartheta(1-c_1)\delta}{L\kappa_d^2}\right\}>0,
\ \
f(\x_k)-f(\x_{k+1})
\geq c_1\delta\underline\alpha_{\mathrm F}\norm{\g_k}^2.
\end{equation}
For a cached Armijo point, either \(\alpha\geq\alpha_{\min}\) or the
secant-type certificate, together with Lipschitz continuity, yields
\(\alpha\geq2\vartheta(1-c_1)\delta/(L\kappa_d^2)\).
Moreover, an Armijo-violating step satisfies
\(\alpha_k^{\mathrm{bad}}>2(1-c_1)\delta/(L\kappa_d^2)\).
Since \(\alpha_k^{(0)}\geq\alpha_{\min}\), reducing the initial seed to
\eqref{eq:ls-m2-reduced-seed} preserves the lower bound \(\underline{\alpha}_{\mathrm F}\). Hence, the fallback-decrease property holds for all accepted fallback steps in exact arithmetic. Under the same assumptions and fixed parameters, the number of backtracking trials and the associated evaluation counts are uniformly bounded.

Importantly, the neural correction changes only the initial trial supplied
to the verified fallback procedure. Neither the corrected proposal nor the
original value \(a_k^{\min}\) is automatically accepted. A line-search
failure does not discard the selected \texttt{CG} direction or
automatically use the \texttt{SD} direction defined in \eqref{eq:steepest-descent-direction}. A fallback iterate has the form \eqref{eq:iteration} and is generated only after either the returned line-search point has been certified or the Armijo fallback verification described above has succeeded.

\paragraph{{\bf Optional emergency continuation}}
The implementation contains an optional emergency continuation that is
considered only when the verified Armijo fallback cannot produce an accepted
step. This branch is heuristic and is not part of the convergence or
evaluation-complexity guarantees established below.

Let \(\x^{\mathrm{best}}\) denote the best accepted point stored so far and
let \(f^{\mathrm{best}}\) be its objective value. The emergency branch
considers the finite cached line-search point, the existing
neural-adjusted fallback seed, and the original deterministic fallback
seed. Let \(\epsilon_{\mathrm H}>0\) denote the scaled objective allowance used
by the optional emergency continuation. A candidate is eligible only when its point, objective value, and
gradient are finite and
\[
f(\x^{\mathrm{trial}})
\leq
f^{\mathrm{best}}
+
\epsilon_{\mathrm H}
\max\left\{
1,
|f^{\mathrm{best}}|
\right\}.
\]

Among the eligible candidates, the point with the smallest objective value
is selected. If two candidates have the same objective value, the one
closest to the stored best point is preferred. If no distinct candidate is
usable, the stored best point is retained.

This emergency continuation may increase the objective relative to the
current iterate, may leave the iterate unchanged, or may return a point that
does not lie on the current search ray. It therefore carries no Armijo or
Wolfe acceptance claim and supplies no fallback-step-\texttt{ELM} training
sample. The emergency branch is an optional implementation mechanism and lies
outside the convergence and evaluation-complexity guarantees.

Under the smoothness and uniform-descent assumptions used in the theoretical
analysis, finite Armijo backtracking succeeds in exact arithmetic.
Consequently, the emergency branch is inactive in that model.

\section{Training Targets and Model Fitting}
\label{sec:training}

In this section, we describe the online training strategy employed for the two \texttt{ELM} models incorporated into the learning-assisted framework. In particular, we specify the construction of their training targets, the corresponding replay archives, and the model-fitting procedures. The coefficient and fallback-step \texttt{ELM}s are trained independently from information generated during the optimization process, without requiring any offline data or additional function or gradient evaluations.

\subsection{Online Coefficient Targets}

The coefficient \texttt{ELM} is trained online using coefficients generated by valid classical \texttt{CG} iterations. No candidate-based teacher search or externally generated target is employed. Recall the analytical coefficient \(\beta_k^{\mathrm{raw}}\) from \eqref{eq:raw-beta-method} and the corresponding classical direction \(\dvec_k^{\mathrm{class}}\) from \eqref{eq:classical-cg-direction}. A coefficient is admitted to the replay archive only when the analytical coefficient is finite and nonzero, the corresponding raw classical direction is finite and strictly descending, and no method-specific restart condition is active. The training sample is stored before application of the final common direction guard. Consequently, the raw classical training direction need not already satisfy the \texttt{DSC} in
\eqref{eq:ls-m3-direction-bounds}. The learning target is the analytical coefficient itself, rather than its subsequently safeguarded value. Under these conditions, the feature-target pair \(\left(\bu_k^{\beta},r_k^{\beta}\right)\) is inserted into the coefficient replay archive. Here, \(\bu_k^{\beta}\in\R^{12}\) is the coefficient feature vector given by \eqref{eq:nn-input}, while \(r_k^{\beta}\) is the encoded classical coefficient defined by
\begin{equation}
	\label{eq:beta-elm-target}
	r_k^{\beta}:=\operatorname{atanh}\left(\operatorname{clip}\left({\beta_k^{\mathrm{raw}}}/{\beta_{\max}},-1+\varepsilon_{\beta},1-\varepsilon_{\beta}\right)\right),
\end{equation}
where \(\beta_{\max}>0\) and \(0<\varepsilon_{\beta}<1\) are tuning parameters, and  $\operatorname{clip}$ is from \eqref{eq:clip-definition}. Hence, the target is a bounded transformation of the analytical coefficient itself, rather than a residual with respect to another \texttt{CG} formula.

\subsection{Coefficient \texttt{ELM} Fitting}

The first \texttt{ELM} is associated with the nonlinear \texttt{CG} coefficient. Its input is the dimension-independent feature vector \(\bu_k^{\beta}\in\R^{12}\) defined by  \eqref{eq:nn-input}. Let \(p_{\beta}\) denote the input dimension and \(h_{\beta}\) the number
of neurons in the single hidden layer. The random input-weight matrix and hidden-layer bias vector are initialized once using a fixed implementation seed and remain unchanged throughout the optimization run. Consequently, each online model update modifies only the output-weight vector.

For a valid classical coefficient \(\beta_k^{\mathrm{raw}}\), the supervised target is formed by first scaling the coefficient by \(\beta_{\max}\) and then applying the inverse hyperbolic tangent, as described in \eqref{eq:beta-elm-target}. Equation \eqref{eq:beta-elm-prediction} yields an absolute \texttt{CG} coefficient rather than a residual correction to the analytical formula. The analytical coefficient remains the primary coefficient. The \texttt{ELM} output is considered only by the restart-recovery or persistent-stagnation mechanisms described in
Section~\ref{sec:nn-parameters}, and any coefficient proposed by the \texttt{ELM} remains subject to the deterministic descent safeguards. Moreover, the coefficient \texttt{ELM} is trained online using valid classical coefficients.

Whenever the coefficient archive satisfies the relevant update conditions, the current archive is used to refit the \texttt{ELM} output layer. The archived feature vectors are standardized and clipped according to \eqref{eq:elm-standardization}--\eqref{eq:elm-feature-clipping}. The fixed hidden-layer representation is then constructed using \eqref{eq:elm-hidden-layer}, and the corresponding hidden matrix is formed as in \eqref{eq:elm-hidden-matrix}. The output-weight vector is obtained by solving the ridge-regression system \eqref{eq:elm-ridge-training}, using the regularization matrix \eqref{eq:elm-ridge-matrix}. Thus, the online fitting procedure updates only the output weights associated with the current coefficient-training archive, while the randomly initialized input weights and hidden-layer biases remain fixed throughout the solver run. This corresponds to standard ridge-regularized least-squares fitting; no Huber loss, sample weighting, candidate merit function, or teacher-generated coefficient is used.

\subsection{Training of the Line-Search Fallback \texttt{ELM}}

A second \texttt{ELM} is maintained independently for the fallback step used after \texttt{cvsrch} failure. Its training samples are obtained from accepted points returned by \texttt{cvsrch}. These include ordinary certified Wolfe points and returned points accepted as certified cached Armijo fallbacks. Fresh Armijo-backtracking points and points selected by the optional emergency continuation are not inserted into this replay archive. Consequently, construction of the fallback-step training targets requires no additional objective-function or gradient evaluations beyond quantities already generated at the accepted returned line-search point.

The input to this model is a ten-dimensional feature vector. Let \(p_{\alpha}\) denote the input dimension of the fallback-step \texttt{ELM}, and let \(h_{\alpha}\) denote the number of hidden units in its single hidden layer. Also, let \(\dvec_k\) denote the direction submitted to \texttt{cvsrch}. For a finite, nonzero direction, define \(\bar{\dvec}_k:=\dvec_k/\norm{\dvec_k}\), and let \(a_k^{\min}\) denote the step obtained from the original fallback formula \eqref{eq:original-amin}. If \(\x_{k+1}^{\mathrm{LS}}\) is an accepted point returned by \texttt{cvsrch}, the corresponding accepted distance along \(\bar{\dvec}_k\) is
\begin{equation}
	\label{eq:accepted-normalized-step}
	\alpha_k^{\mathrm{acc}}:=\left(\x_{k+1}^{\mathrm{LS}}-\x_k\right)^{\mathsf T}\bar{\dvec}_k.
\end{equation}
A training sample is retained only when \(a_k^{\min}\) and \(\alpha_k^{\mathrm{acc}}\) are finite and strictly positive.

% Accordingly, the implementation uses
% \begin{equation}
% 	\label{eq:step-elm-hidden-units}
% 	p_{\alpha}=10,\qquad m_{\alpha}=16.
% \end{equation}

Unlike the coefficient \texttt{ELM}, the output of this model is interpreted as a logarithmic step length correction. The corresponding training target is defined by
\begin{equation}
	\label{eq:step-elm-target}
	r_k^{\alpha}:=\operatorname{clip}\left(\log({\alpha_k^{\mathrm{acc}}}/{a_k^{\min}}),-\log C_{\mathrm{tar}},\log C_{\mathrm{tar}}\right),
\end{equation}
where \(C_{\mathrm{tar}}>1\) is a prescribed constant and $\operatorname{clip}$ is from \eqref{eq:clip-definition}.
% \begin{equation}
% 	\label{eq:step-target-correction}
% 	C_{\mathrm{tar}}=10^8.
% \end{equation}
Using the regularized cosine \(C(\cdot,\cdot)\) from \eqref{eq:feature-cosine} and the signed logarithm \(L(\cdot)\) from \eqref{eq:signed-log}, the ten-dimensional feature vector is now defined by
\begin{equation}
	\label{eq:step-feature-vector}
	\begin{aligned}
	\bu_k^\alpha:=\Big[&
	L(f_k),\,
	\log\!\left(1+\norm{\g_k}\right),\,
	C(\g_k,\bar{\dvec}_k),\\
	&
	L\!\left(\g_k^{\mathsf T}\bar{\dvec}_k\right),\,
	\log\!\left(1+\norm{\x_k}_{\infty}\right),\,
	L(\xi_k),\\
	&
	\log\!\left(\max\{a_k^{\min},\mathrm{realmin}\}\right),\,
	L\!\left(\beta_{k}^{\mathrm{used}}\right),\,
	P(R_k^f),\,P(R_k^g)
	\Big]^{\mathsf T},
	\end{aligned}
\end{equation}
where 
\begin{equation}
	\label{eq:step-x-over-d}
	\xi_k:=
	\begin{cases}
		\displaystyle\min_{\left\{i:\ \bar d_{k,i}\neq0\right\}}\left|x_{k,i}/\bar d_{k,i}\right|, &\qquad \text{if at least one }\bar d_{k,i}\neq0,\\[1ex]
		0, &\qquad \text{otherwise},
	\end{cases}
\end{equation}
and
\begin{equation}
	\label{eq:progress-step-transform}
	P(v):=
	\begin{cases}
		0, &\qquad v\text{ is nonfinite or }v<0,\\
		\log(1+v), &\qquad \text{otherwise}.
	\end{cases}
\end{equation}
The quantities \(R_k^f\) and \(R_k^g\) in \eqref{eq:step-feature-vector}
are the stagnation measures defined in \eqref{eq:relative-f-variation}
and \eqref{eq:relative-g-variation}, respectively. Although these quantities are not components of the coefficient-\texttt{ELM} feature vector, they provide the fallback-step model with additional information on recent optimization progress. Any nonfinite component of \(\bu_k^\alpha\) is replaced by zero, after which the feature vector is clipped using the same feature bound introduced in \eqref{eq:elm-feature-clipping}. Model fitting then follows the generic \texttt{ELM} standardization and ridge-regression procedure \eqref{eq:elm-standardization}--\eqref{eq:elm-ridge-training}.

% The fallback-step \texttt{ELM} uses an archive of at most \(250\) samples, begins fitting after \(8\) successful line search samples, is refitted after every \(4\) new samples, and uses the ridge parameter \(\lambda_{\alpha}=10^{-4}\). Its hidden weights and biases are initialized independently of those of the coefficient \texttt{ELM}. With the implementation seeds used in the solver, the coefficient-\texttt{ELM} seed is \(73\), whereas the fallback-step-\texttt{ELM} seed is \(73+1009=1082\).

At prediction time, the fallback-step-\texttt{ELM} produces
\(\widehat r_k^\alpha\), whose multiplicative interpretation
\(C_k^\alpha\) is defined in \eqref{eq:step-elm-correction}. The fallback-step model applies the standardized-feature {\tt OOD} screening introduced in \eqref{eq:ood-test}, using the standardization statistics and {\tt OOD} scale computed from its own replay archive. The predicted logarithmic correction is restricted to the same interval used for the target in \eqref{eq:step-elm-target}, while the additional safeguards in \eqref{eq:amin-correction-bounds}
restrict the multiplicative correction. Thus, the fallback-step \texttt{ELM} can enlarge the original fallback step but cannot reduce it.

Accordingly, the two \texttt{ELM} models are trained online from distinct information sources. Valid classical \texttt{CG} coefficients provide the training samples for the coefficient \texttt{ELM}, whereas accepted points returned by \texttt{cvsrch}, including certified cached Armijo points, provide the training samples for the fallback-step \texttt{ELM}. Fresh fallback-backtracking points and emergency-continuation points are excluded from the latter archive. The two models have separate feature vectors, replay archives, standardization statistics, hidden-layer parameters, and output weights. Neither model requires an additional objective-function or gradient evaluation solely for training.

\section{Computational Cost of the Learning-Assisted Mechanism}
\label{sec:computational-cost}

The proposed learning-assisted variants introduce additional algebraic operations relative to their classical \texttt{CG} counterparts. These operations arise from constructing the neural features, evaluating the two \texttt{ELM} models, updating their replay archives, fitting their output coefficients, performing the \texttt{OOD} tests, and applying the deterministic safeguards. In the present implementation, these neural operations do not themselves require additional objective-function or gradient evaluations. The neural models modify either the \texttt{CG} coefficient or, after a line-search failure, the fallback step that would otherwise be used by the deterministic method.

The implementation employs two independent online \texttt{ELM} models. The
first predicts a safeguarded \texttt{CG} coefficient and uses a
\(12\)-dimensional feature vector. The second predicts a correction of the
fallback step after a line-search failure and uses a \(10\)-dimensional
feature vector. Both models have fixed randomly generated hidden-layer
parameters and determine only their output coefficients through regularized
least squares. In what follows, \(\mathcal{C}\) denotes computational cost;
the subscript specifies the operation being performed, while the superscript
\(\beta\) or \(\alpha\) identifies, respectively, the coefficient-\texttt{ELM}
or the fallback-step-\texttt{ELM}.

\subsection{Feature Construction Cost}

Let \(p_\beta=12\) denote the input dimension of the coefficient-\texttt{ELM}.
Its feature vector contains a fixed number of norms, inner products, cosine
quantities, and scalar transformations involving the current and previous
optimization quantities. Let \(\mathcal{C}_{\mathrm{feat}}^\beta\) denote
the computational cost of constructing this feature vector. Since only a
fixed number of operations on vectors in \(\R^n\) is required, $\mathcal{C}_{\mathrm{feat}}^\beta=\mathcal{O}(n)$. 

The fallback-step model has input dimension \(p_\alpha=10\). Its feature
construction contains norms, inner products, an infinity norm, and a
coordinatewise calculation involving the current point and the normalized
search direction. Let \(\mathcal{C}_{\mathrm{feat}}^\alpha\) denote the
corresponding feature-construction cost. Hence, $\mathcal{C}_{\mathrm{feat}}^\alpha=\mathcal{O}(n)$. Both \(p_\beta\) and \(p_\alpha\) are independent of the optimization
dimension \(n\). Therefore, increasing \(n\) does not increase the input
dimension of either neural model.

\subsection{Cost of One \texttt{ELM} Prediction}

Consider the generic \texttt{ELM} introduced in
\Cref{subsec:nn-preliminaries}, with input dimension \(p\) and \(h\) hidden
units. Standardization, clipping, and \texttt{OOD} screening require
\(\mathcal{O}(p)\) operations. Computing the hidden representation requires
the multiplication of a \(p\)-dimensional input by a \(p\times h\)
hidden-weight matrix and therefore costs \(\mathcal{O}(ph)\). The scalar
output evaluation requires only \(\mathcal{O}(h)\) additional operations.
Consequently,
\begin{equation}
\label{eq:elm-prediction-cost}
\mathcal{C}_{\mathrm{pred}}(p,h)=\mathcal{O}(ph).
\end{equation}
For the coefficient-\texttt{ELM}, \(p_\beta=12\), while \(h_\beta\)
denotes its fixed hidden-layer size. Thus,
\(\mathcal{C}_{\mathrm{pred}}^\beta=\mathcal{O}(p_\beta h_\beta)\).
For the fallback-step-\texttt{ELM}, \(p_\alpha=10\), while \(h_\alpha\)
denotes its fixed hidden-layer size, giving
\(\mathcal{C}_{\mathrm{pred}}^\alpha=\mathcal{O}(p_\alpha h_\alpha)\). Since \(p_\beta,h_\beta,p_\alpha\), and \(h_\alpha\) are fixed independently
of \(n\), the direct cost of either neural prediction is
\(\mathcal{O}(1)\) with respect to the optimization dimension.
% \begin{equation}
% \label{eq:beta-prediction-cost}
% \mathcal{C}_{\mathrm{pred}}^\beta=\mathcal{O}(p_\beta h_\beta).
% \end{equation}
% \begin{equation}
% \label{eq:alpha-prediction-cost}
% \mathcal{C}_{\mathrm{pred}}^\alpha=\mathcal{O}(p_\alpha h_\alpha).
% \end{equation}

The subsequent construction and checking of the associated
\texttt{CG} direction require vector operations in \(\R^n\). Since the
number of descent-repair operations is bounded by a fixed implementation
constant, the complete algebraic cost associated with producing and
safeguarding a neural direction remains $\mathcal{O}(n)$.
% \begin{equation}
% \label{eq:nn-safeguard-cost}
% \mathcal{O}(n).
% \end{equation}

\subsection{Cost of One \texttt{ELM} Fitting Event}

Suppose that an \texttt{ELM} replay archive currently contains \(N\) samples,
with \(p\) input features and \(h\) hidden units. Computing the feature means
and standard deviations requires \(\mathcal{O}(Np)\) operations and forming the
hidden-layer matrix requires
\(\mathcal{O}(Nph)\) operations. The implementation determines the output coefficients by solving the
regularized least-squares system associated with the hidden-layer matrix.
Formation of its normal matrix requires \(\mathcal{O}(Nh^2)\) operations,
while solution of the resulting \((h+1)\times(h+1)\) dense linear system
requires \(\mathcal{O}(h^3)\) operations.

After fitting, the standardized row magnitudes are sorted to determine the
reference scale used by the \texttt{OOD} test. This operation contributes
\(\mathcal{O}(N\log N)\). Therefore, one complete model-fitting event 
costs as
\begin{equation}
\label{eq:elm-fitting-cost}
\mathcal{C}_{\mathrm{fit}}(N,p,h)=\mathcal{O}(Nph+Nh^2+h^3+N\log N).
\end{equation}
Unlike a hyperparameter-tuning procedure, the present implementation does
not perform a grid search over hidden-layer sizes or ridge parameters during
a fitting event. The hidden-layer size and ridge parameter are fixed for
each of the two models.

\subsection{Cost of the Coefficient-\texttt{ELM}}

Let \(N_\beta\) denote the current number of samples in the coefficient
archive, whose fixed capacity is denoted by \(N_\beta^{\max}\). Substituting these quantities into \eqref{eq:elm-fitting-cost} gives
\begin{equation}
\label{eq:beta-elm-fitting-cost}
\mathcal{C}_{\mathrm{fit}}^\beta=
\mathcal{O}(N_\beta p_\beta h_\beta+N_\beta h_\beta^2+h_\beta^3+N_\beta\log N_\beta).
\end{equation}
The minimum number of samples required for the first fit, the retraining
interval, and the archive capacity are fixed independently of the
optimization dimension. Since \(N_\beta,p_\beta\), and \(h_\beta\) are all bounded independently of
\(n\), the cost of each coefficient-\texttt{ELM} update is bounded with
respect to the optimization dimension.

\subsection{Cost of the Fallback-Step-\texttt{ELM}}

Let \(N_\alpha\) denote the number of samples in the fallback-step archive,
whose fixed capacity is denoted by \(N_\alpha^{\max}\). Thus,
\begin{equation}
\label{eq:alpha-elm-fitting-cost}
\mathcal{C}_{\mathrm{fit}}^\alpha=
\mathcal{O}(N_\alpha p_\alpha h_\alpha+N_\alpha h_\alpha^2+h_\alpha^3+N_\alpha\log N_\alpha).
\end{equation}

As for the coefficient model, the minimum fitting size, retraining interval,
hidden-layer size, and ridge parameter remain fixed throughout the
optimization. The fallback-step training samples are obtained from successful
\texttt{cvsrch} iterations. Their construction uses values already produced
by the optimization method and therefore requires no additional objective-function or
gradient evaluation.

\subsection{Total Algebraic Learning Cost}

Let \(K\) denote the number of \texttt{CG} iterations,
\(K_\beta\) and \(K_\alpha\) the numbers of coefficient and fallback-step
predictions, respectively, and \(U_\beta\) and \(U_\alpha\) the corresponding
numbers of model updates. The total direct algebraic contribution of the
learning mechanisms can be written as
\begin{equation}
\label{eq:total-learning-algebraic-cost}
\mathcal{C}_{\mathrm{LA}}^{\mathrm{alg}}=
\mathcal{O}(Kn)+\mathcal{O}(K_\beta p_\beta h_\beta)+
\mathcal{O}(K_\alpha p_\alpha h_\alpha)+
U_\beta\mathcal{C}_{\mathrm{fit}}^\beta+
U_\alpha\mathcal{C}_{\mathrm{fit}}^\alpha.
\end{equation}

At most one new sample can be added to either archive during one optimization
iteration, and each model is updated only after a fixed number of newly
collected samples. Hence, $U_\beta=\mathcal{O}(K)$, and $U_\alpha=\mathcal{O}(K)$.
Since the archive capacities, feature dimensions, and hidden-layer sizes are
fixed independently of \(n\), the model-update and prediction terms in
\eqref{eq:total-learning-algebraic-cost} are bounded independently of \(n\)
per iteration. Therefore,
\begin{equation}
\label{eq:learning-algebraic-order}
\mathcal{C}_{\mathrm{LA}}^{\mathrm{alg}}=\mathcal{O}(Kn)
\end{equation}
with respect to the optimization dimension. Thus, the learning mechanisms do not alter the dimension-dependent
arithmetic order of the underlying nonlinear \texttt{CG} iteration.

\subsection{Storage Requirements}

The classical nonlinear \texttt{CG} method stores a fixed number of vectors
in \(\R^n\) and therefore requires \(\mathcal{O}(n)\) storage. The
learning-assisted implementation additionally stores two bounded archives,
two fixed hidden-layer models, normalization quantities, and short scalar
histories. The coefficient archive requires
\(\mathcal{O}(N_\beta^{\max}p_\beta)\) storage and the fallback-step archive
requires
\(\mathcal{O}(N_\alpha^{\max}p_\alpha)\). The two hidden-weight matrices require
\(\mathcal{O}(p_\beta h_\beta)\) and
\(\mathcal{O}(p_\alpha h_\alpha)\) storage, respectively. All of these
quantities are independent of \(n\). Consequently, the total storage requirement remains
\begin{equation}
\label{eq:nncg-storage-complexity}
\mathcal{O}(n).
\end{equation}

\begin{proposition}
\label{prop:learning-cost}
Suppose that the neural feature dimensions, hidden-layer sizes, archive
capacities, retraining intervals, and maximum number of descent-repair steps
are fixed independently of \(n\). Then the direct algebraic overhead
introduced by the learning-assisted mechanisms is
\(\mathcal{O}(n)\) per nonlinear \texttt{CG} iteration, and the total
storage complexity remains \(\mathcal{O}(n)\).
\end{proposition}

\begin{proof}
The construction of the two feature vectors and the deterministic direction
safeguards require only a fixed number of vector operations in \(\R^n\) and
therefore cost \(\mathcal{O}(n)\). By
\eqref{eq:elm-prediction-cost}, neural prediction has fixed cost with respect
to \(n\). By \eqref{eq:beta-elm-fitting-cost} and
\eqref{eq:alpha-elm-fitting-cost}, each fitting event is also bounded with
respect to \(n\), because both archives and both neural models have fixed
sizes. Hence the additional dimension-dependent work per iteration is
\(\mathcal{O}(n)\). The two replay archives and all neural-model parameters have sizes independent
of \(n\). Since the underlying \texttt{CG} method already stores a fixed
number of \(n\)-dimensional vectors, adding the neural data structures leaves
the total storage order equal to \(\mathcal{O}(n)\).
\end{proof}

\begin{remark}
\label{rem:evaluation-accounting}
The algebraic cost described above should be distinguished from function and
gradient evaluations. Neither fitting nor evaluating either \texttt{ELM}
requires an additional call to the objective function or its gradient. The
coefficient model changes the coefficient used to form the search direction,
while the fallback-step model changes the fallback step before the fallback
point is evaluated. Their effect on the evaluation counts is therefore
indirect, through the trajectory, the behavior of the line search, and the
fallback step selected by the solver.
\end{remark}

\subsection{Convergence and Evaluation Complexity of the Learning-Assisted
\texttt{CG} Method}
\label{subsec:nn-evaluation-complexity}

The analysis below concerns the verified learning-assisted algorithm in
exact arithmetic and in the absence of external stopping limits. Under the
smoothness and uniform-descent assumptions stated below, finite Armijo
backtracking succeeds. Consequently, the optional emergency continuation is
inactive in the theoretical model. The results therefore apply to
iterations generated by the verified Wolfe, certified cached-Armijo, or
fresh Armijo-fallback branches and do not assert the same whole-run
guarantees for an actual trajectory containing an unverified emergency
continuation.

The convergence analysis depends on the final safeguarded direction and
accepted step rather than on the source of the coefficient. The coefficient
may originate from a classical formula, the coefficient-\texttt{ELM}, the
recent coefficient history, or the deterministic safe bootstrap.

Before \texttt{cvsrch} is called, every direction produced by a
learning-assisted variant is subjected to the common final direction guard.
Hence the \texttt{SDC} and \texttt{DSC} used below are enforced
algorithmically. Theorem~\ref{thm:ls-m3-simultaneous}
establishes the corresponding simultaneous bounds in exact arithmetic,
while the final implementation check enforces their finite-precision
counterparts.

Likewise, an ordinary successful \texttt{cvsrch} iteration is accepted only
after the corresponding line-search conditions have been certified. Returned
points accepted as cached Armijo fallbacks and points produced by fresh
fallback backtracking satisfy the Armijo condition by their respective
acceptance tests. The objective function is assumed to satisfy the following standard
condition \cite{nocedal1999numerical}.

\begin{assumption}
\label{ass:bounded-below-lipschitz}
The objective function \(f\) is bounded below on a convex set containing the
iterates and all line-search and fallback trial points; that is, there
exists \(f_{\inf}>-\infty\) such that $f(\x)\geq f_{\inf}$. Moreover, the gradient of \(f\) is Lipschitz continuous on this set; that
is, there exists \(L>0\) such that
\[
\norm{\g(\x)-\g(\bz)}
\leq
L\norm{\x-\bz}
\]
for all relevant \(\x,\bz\) in this set.
\end{assumption}

\begin{lemma}
\label{lem:wolfe-step-lower-bound}
Suppose that Assumption~\ref{ass:bounded-below-lipschitz} holds.
Then every successful Wolfe step satisfies
\begin{equation}
\label{eq:alpha-lower-bound}
\alpha_k
\geq
\underline{\alpha}
:=
\frac{(1-c_2)\delta}{L\kappa_d^2}
>0.
\end{equation}
\end{lemma}

\begin{proof}
From the Wolfe curvature condition \eqref{eq:wolfe-curvature},
\[
\g_{k+1}^{\mathsf T}\dvec_k
-
\g_k^{\mathsf T}\dvec_k
\geq
(1-c_2)
\left(
-\g_k^{\mathsf T}\dvec_k
\right).
\]
By Lipschitz continuity of the gradient and \eqref{eq:iteration},
\[
(\g_{k+1}-\g_k)^{\mathsf T}\dvec_k
\leq
\norm{\g_{k+1}-\g_k}\norm{\dvec_k}
\leq
L\alpha_k\norm{\dvec_k}^{2}.
\]
Therefore,
\[
\alpha_k
\geq
\frac{
(1-c_2)
(-\g_k^{\mathsf T}\dvec_k)
}{
L\norm{\dvec_k}^{2}
}.
\]
Applying the \texttt{SDC} and \texttt{DSC} in
\eqref{eq:ls-m3-direction-bounds} leads to \eqref{eq:alpha-lower-bound}.
\end{proof}

\begin{lemma}
\label{lem:uniform-wolfe-decrease}
Suppose that Assumption~\ref{ass:bounded-below-lipschitz} holds, and define
\begin{equation}
\label{eq:complexity-mu-wolfe}
\mu_{\mathrm W}
:=
\frac{c_1(1-c_2)\delta^2}
{L\kappa_d^2}
>0.
\end{equation}
Then every successful Wolfe iteration satisfies
\begin{equation}
\label{eq:uniform-wolfe-decrease}
f_k-f_{k+1}
\geq
\mu_{\mathrm W}\norm{\g_k}^{2}.
\end{equation}
\end{lemma}

\begin{proof}
By the Armijo condition \eqref{eq:wolfe-armijo}, $f_k-f_{k+1}
\geq
-c_1\alpha_k
\g_k^{\mathsf T}\dvec_k$. Using the \texttt{SDC} in
\eqref{eq:ls-m3-direction-bounds} and
Lemma~\ref{lem:wolfe-step-lower-bound}, we obtain
\[
f_k-f_{k+1}
\geq
c_1\delta\alpha_k\norm{\g_k}^{2}
\geq
\mu_{\mathrm W}\norm{\g_k}^{2}.
\]
\end{proof}

The fallback branch satisfies an analogous uniform decrease property without
requiring an additional convergence assumption. Indeed, the Armijo
backtracking safeguard \eqref{eq:ls-m2-uniform-decrease} already guarantees
the required decrease for every accepted fallback step. Define $\mu_{\mathrm F}
:=
c_1\delta\underline{\alpha}_{\mathrm F}
>0$. Then, by \eqref{eq:ls-m2-uniform-decrease}, every accepted fallback
iteration satisfies
\begin{equation}
\label{eq:fallback-mu-decrease}
f_k-f_{k+1}
\geq
\mu_{\mathrm F}\norm{\g_k}^{2}.
\end{equation}
Thus, the fallback decrease used in the convergence analysis is a direct
consequence of \eqref{eq:ls-m2-uniform-decrease}, rather than an independent
assumption.

\begin{theorem}
\label{thm:nn-global-convergence}
Suppose that Assumption~\ref{ass:bounded-below-lipschitz} holds.
Then either the method terminates at a stationary point in finitely many
iterations, or the generated sequence satisfies
\begin{equation}
\label{eq:gradient-square-summability}
\sum_{k=0}^{\infty}
\norm{\g_k}^{2}
<\infty,
\end{equation}
and consequently
\begin{equation}
\label{eq:gradient-convergence}
\lim_{k\rightarrow\infty}
\norm{\g_k}
=
0.
\end{equation}
\end{theorem}

\begin{proof}
If the method terminates with \(\g_k=0\), the first conclusion holds.
Otherwise, consider the infinite sequence of accepted iterations.

By Lemma~\ref{lem:uniform-wolfe-decrease}, every successful Wolfe
iteration satisfies \eqref{eq:uniform-wolfe-decrease}. By the fallback safeguard \eqref{eq:ls-m2-uniform-decrease} and the
definition of \(\mu_{\mathrm F}\), every accepted fallback iteration
satisfies \eqref{eq:fallback-mu-decrease}. Therefore, every accepted iteration satisfies
\[
f_k-f_{k+1}
\geq
\bar{\mu}\norm{\g_k}^{2} \ \ \text{with $\bar{\mu}
	:=
	\min\{\mu_{\mathrm W},\mu_{\mathrm F}\}
	>0.$}
\]
Summing from \(k=0\) to \(K-1\) yields
\[
\bar{\mu}
\sum_{k=0}^{K-1}
\norm{\g_k}^{2}
\leq
f_0-f_K
\leq
f_0-f_{\inf}.
\]
Letting \(K\rightarrow\infty\) proves
\eqref{eq:gradient-square-summability}. Since all summands are
nonnegative, it follows that
\(\norm{\g_k}\rightarrow0\).
\end{proof}

\begin{theorem}
\label{thm:nn-iteration-complexity}
For \(\varepsilon>0\), let \(K_\varepsilon\) denote the first iteration
index for which
\[
\norm{\g_{K_\varepsilon}}_\infty
\leq
\varepsilon.
\]
Under the conditions of
Theorem~\ref{thm:nn-global-convergence},
\begin{equation}
\label{eq:iteration-complexity-bound}
K_\varepsilon
\leq
\left\lceil
{
(f_0-f_{\inf})
}/{
(\bar{\mu}\varepsilon^2)
}
\right\rceil.
\end{equation}
Consequently,
\begin{equation}
\label{eq:iteration-complexity-order}
K_\varepsilon
=
\mathcal{O}(\varepsilon^{-2}).
\end{equation}
\end{theorem}

\begin{proof}
Suppose that $\norm{\g_k}_\infty>\varepsilon$ for all
$k=0,\ldots,K-1$. Since
\(\norm{\g_k}\geq\norm{\g_k}_\infty\), $\norm{\g_k}^{2}
>
\varepsilon^{2}$. The uniform decrease established in
Theorem~\ref{thm:nn-global-convergence} therefore gives
$f_k-f_{k+1}
>
\bar{\mu}\varepsilon^{2}$. Summing over the first \(K\) iterations yields $f_0-f_K
>
K\bar{\mu}\varepsilon^{2}$. Since \(f_K\geq f_{\inf}\), $K
<{
(f_0-f_{\inf})
}/{
(\bar{\mu}\varepsilon^{2})
}$, which proves \eqref{eq:iteration-complexity-bound}.
\end{proof}

We next translate the iteration bound into evaluation complexity. In the
present implementation, each trial performed internally by
\texttt{cvsrch} evaluates the objective function and gradient jointly.
Let \(\ell_{\max}\) denote the maximum number of such joint evaluations
in one call to \texttt{cvsrch}.

The fallback branch requires separate accounting because its Armijo
verification may evaluate several trial points before a fallback step is
accepted. When separate objective and gradient routines are available,
rejected fallback trials require only objective-function evaluations and
the gradient is evaluated only after Armijo acceptance. When only a joint
function-gradient routine is available, every fallback trial evaluates
both the objective function and the gradient. The corresponding evaluation
bounds therefore depend on which of these two access models is available.

\begin{lemma}
\label{lem:nn-no-extra-evaluations}
Neither the coefficient-\texttt{ELM} nor the fallback-step-\texttt{ELM}
requires an additional objective-function or gradient evaluation solely for
feature construction, prediction, archive updating, or model fitting.
\end{lemma}

\begin{proof}
The coefficient-\texttt{ELM} uses quantities already available from the
current and previous optimization iterations. Its prediction only modifies
the coefficient used to form the search direction. The fallback-step model
is trained using accepted points already returned by \texttt{cvsrch},
including certified cached Armijo points, and only modifies the initial
fallback proposal. Hence neither model requires an objective-function or
gradient evaluation solely for learning.

The function evaluations used by fallback Armijo backtracking belong to
the deterministic fallback acceptance safeguard and are therefore accounted
for separately below.
\end{proof}

We next bound the number of fallback trial evaluations. From the fallback
construction, the initial trial step satisfies $\alpha_k^{(0)}
\leq
\alpha_{\max}$, and each rejected fallback trial is contracted by the fixed factor
\(\vartheta\in(0,1)\). By the descent-lemma argument leading to
\eqref{eq:ls-m2-uniform-decrease}, the Armijo condition is guaranteed once
\[
\alpha
\leq{
2(1-c_1)\delta
}/{
(L\kappa_d^2)
}.
\]
Consequently, the number of function evaluations performed by one fallback
backtracking sequence is bounded by the finite constant
\begin{equation}
\label{eq:fallback-max-trials}
\ell_{\mathrm F}
:=
1+
\max\left\{
0,\,
\left\lceil
{
\log\!\left(
\dfrac{
\alpha_{\max}L\kappa_d^2
}{
2(1-c_1)\delta
}
\right)
}\Big/{
\log(\vartheta^{-1})
}
\right\rceil
\right\}.
\end{equation}
This bound includes the accepted fallback trial.

Let \(N_{\mathrm{fb}}(K)\) denote the number of iterations among the first
\(K\) iterations for which the fallback branch is entered and an accepted
fallback point is generated. Let
\({\tt nf}_\varepsilon\) and
\({\tt ng}_\varepsilon\) denote, respectively, the numbers of objective
function and gradient evaluations required to obtain
\(\norm{\g_k}_\infty\leq\varepsilon\), and define
\begin{equation}
\label{eq:nf2g}
{\tt nf2g}_\varepsilon
:=
{\tt nf}_\varepsilon
+
2\,{\tt ng}_\varepsilon.
\end{equation}

\begin{theorem}
\label{thm:nf-ng-complexity}
Under the conditions of
Theorem~\ref{thm:nn-iteration-complexity}, the evaluation counts satisfy
\begin{equation}
\label{eq:nf-ng-evaluation-bound}
\begin{aligned}
{\tt nf}_\varepsilon
&\leq
1
+
\ell_{\max}K_\varepsilon
+
\ell_{\mathrm F}
N_{\mathrm{fb}}(K_\varepsilon),
\\
{\tt ng}_\varepsilon
&\leq
1
+
\ell_{\max}K_\varepsilon
+
N_{\mathrm{fb}}(K_\varepsilon).
\end{aligned}
\end{equation}
Since
\(N_{\mathrm{fb}}(K_\varepsilon)\leq K_\varepsilon\),
we obtain
\begin{equation}
\label{eq:nf-ng-simple-bound}
{\tt nf}_\varepsilon
\leq
1+
(\ell_{\max}+\ell_{\mathrm F})
K_\varepsilon,
\ \
{\tt ng}_\varepsilon
\leq
1+
(\ell_{\max}+1)
K_\varepsilon.
\end{equation}
With joint-only fallback access, the gradient bounds in
\eqref{eq:nf-ng-evaluation-bound}--\eqref{eq:nf-ng-simple-bound}
are the corresponding function-evaluation bounds. Consequently,
\begin{equation}
\label{eq:nf-ng-complexity-order}
{\tt nf}_\varepsilon
=
\mathcal{O}(\varepsilon^{-2}),
\qquad
{\tt ng}_\varepsilon
=
\mathcal{O}(\varepsilon^{-2}),
\end{equation}
and
\begin{equation}
\label{eq:nf2g-complexity-order}
{\tt nf2g}_\varepsilon
=
\mathcal{O}(\varepsilon^{-2}).
\end{equation}
\end{theorem}

\begin{proof}
The initial point requires one joint function-gradient evaluation. Each
trial inside \texttt{cvsrch} performs one joint evaluation, and one
\texttt{cvsrch} call uses at most \(\ell_{\max}\) such evaluations. If the fallback branch is entered, at most \(\ell_{\mathrm F}\) additional
objective-function evaluations are required by the Armijo backtracking
sequence. With separate access, the gradient is evaluated only at acceptance; with joint-only access, every trial also evaluates the gradient. Hence the
fallback gradient cost is one with separate access and at most
\(\ell_{\mathrm F}\) with joint-only access, proving the stated bounds. Using \(N_{\mathrm{fb}}(K_\varepsilon)
\leq K_\varepsilon\) gives \eqref{eq:nf-ng-simple-bound}. Since \(\ell_{\max}\) and \(\ell_{\mathrm F}\) are constants independent of
\(\varepsilon\), Theorem~\ref{thm:nn-iteration-complexity} yields
\eqref{eq:nf-ng-complexity-order}. Finally, \eqref{eq:nf2g} implies \eqref{eq:nf2g-complexity-order}.
\end{proof}
Combining the iteration-complexity result with
\eqref{eq:learning-algebraic-order}, we obtain
\begin{equation}
\label{eq:learning-algebraic-epsilon-order}
\mathcal{C}_{\mathrm{LA},\varepsilon}^{\mathrm{alg}}
=
\mathcal{O}(n\varepsilon^{-2}).
\end{equation}

\begin{corollary}
\label{cor:nn-complexity-preservation}
Suppose that Assumption~\ref{ass:bounded-below-lipschitz} holds.
Then the safeguarded learning-assisted method preserves the worst-case
first-order orders
\[
K_\varepsilon
=
\mathcal{O}(\varepsilon^{-2}),
\qquad
{\tt nf}_\varepsilon
=
\mathcal{O}(\varepsilon^{-2}),
\qquad
{\tt ng}_\varepsilon
=
\mathcal{O}(\varepsilon^{-2}), \ \ \text{and} \ \
{\tt nf2g}_\varepsilon
=
\mathcal{O}(\varepsilon^{-2}).
\]
Moreover, the direct learning-related arithmetic work and total storage
satisfy
\[
\mathcal{C}_{\mathrm{LA},\varepsilon}^{\mathrm{alg}}
=
\mathcal{O}(n\varepsilon^{-2}),
\qquad
\mathcal{S}_{\mathrm{CGdNN}}
=
\mathcal{O}(n).
\]
\end{corollary}

\begin{proof}
The convergence and evaluation-complexity analysis depends on the properties
of the final accepted direction and step, not on the source of the proposed
coefficient. By the uniform safeguard
\eqref{eq:ls-m3-direction-bounds} and
Theorem~\ref{thm:ls-m3-simultaneous}, every final direction supplied to the line search satisfies the
\texttt{SDC} and \texttt{DSC} in
\eqref{eq:ls-m3-direction-bounds}, independently of whether the proposed coefficient is analytical,
generated by the coefficient-\texttt{ELM}, obtained from recent coefficient
information, or constructed by the deterministic safe initialization.

Furthermore, successful \texttt{cvsrch} iterations satisfy the Wolfe
conditions, while accepted fallback iterations satisfy the Armijo decrease
guarantee \eqref{eq:ls-m2-uniform-decrease}. Hence the coefficient-generation
mechanism does not alter the deterministic properties used in
Theorems~\ref{thm:nn-global-convergence}--
\ref{thm:nf-ng-complexity}. The stated first-order complexity orders
therefore follow directly. Equation
\eqref{eq:learning-algebraic-epsilon-order} gives the learning-related
arithmetic order, while \eqref{eq:nncg-storage-complexity} gives the storage
result.
\end{proof}

\begin{remark}
\label{rem:nn-complexity-interpretation}
The preceding results are complexity-preservation results and do not imply
that the classical and learning-assisted methods generate identical
trajectories. A learning-assisted direction may alter the number of
iterations, line-search trials, and fallback events observed in practice.
Likewise, the fallback-step-\texttt{ELM} may change the initial step supplied
to the fallback Armijo procedure. Nevertheless, the uniform direction
safeguard \eqref{eq:ls-m3-direction-bounds} and the fallback acceptance
guarantee \eqref{eq:ls-m2-uniform-decrease} ensure that these changes do not
alter the \(\mathcal{O}(\varepsilon^{-2})\) worst-case first-order order.
The learning models themselves require no additional objective-function or
gradient evaluations solely for training or prediction.
\end{remark}

\begin{remark}
\label{rem:implementation-complexity-conditions}
The direction properties used in the preceding proofs are enforced by the
algorithm rather than postulated as additional assumptions. In exact
arithmetic, the final directions satisfy the \texttt{SDC} and \texttt{DSC}
in \eqref{eq:ls-m3-direction-bounds}. In the finite-precision
implementation, every final classical, recovery, or persistent-stagnation
direction is verified against the \texttt{SDC} and the finite-precision
form of the \texttt{DSC} described above before \texttt{cvsrch} is
called. Similarly, a fallback point is accepted only after satisfying the
Armijo test used to establish \eqref{eq:ls-m2-uniform-decrease}. Thus,
apart from Assumption~\ref{ass:bounded-below-lipschitz}, the descent,
direction-scale, Wolfe-success, and fallback-decrease properties invoked
in the exact-arithmetic analysis follow from the deterministic acceptance
safeguards of the method.
\end{remark}
\section{Numerical Results}
\label{sec:experiments}

The numerical study compares the six classical methods with their neural
counterparts on the 480 \texttt{CUTEst} test problems \cite{CUTEst}, whose
dimensions range from 2 to 9000. In addition, we consider a
PDE-constrained conductivity-identification problem to assess the
performance of \texttt{HS} and \texttt{HSNN} when function and gradient
evaluations are computationally expensive. All methods use identical
starting points, stopping criteria, evaluation budgets, and line-search
settings; only the stated learning and safeguarding mechanisms differ.

A run is declared successful when
\(\norm{\g}_{\infty}\leq 10^{-6}\) under the evaluation budget ${\tt nf2g}\leq 500000$, where \({\tt nf2g}\) denotes the weighted evaluation count corresponding
to \eqref{eq:nf2g}, with the stopping-tolerance subscript omitted for
individual numerical runs. No finite time limit is imposed. For each run, we record the numbers of function and gradient evaluations,
\(\texttt{nf}\) and \(\texttt{ng}\), respectively, together with the weighted
evaluation cost \(\texttt{nf2g}:=\texttt{nf}+2\,\texttt{ng}\) and the elapsed
time \(\texttt{msec}\), measured in milliseconds. Lower values are preferable
for all four measures.

\subsection{Implementation Parameters}
\label{subsec:numerical-parameters}

The implementation fixes \(t_{\mathrm{DL}}=0.1\) for
\texttt{DL}/\texttt{DLNN} \cite{DaiLiaoNCG}, retains the analytical
\texttt{HZ} parameter, and uses the \texttt{DY}/\texttt{FR} restart
criterion \eqref{eq:dy-fr-restart}. Table~\ref{tab:cg-wolfe-parameters} in \ref{sec:appendix1} lists the classical, stopping,
fallback-seed, and Wolfe settings. Table~\ref{tab:nn-parameters} in \ref{sec:appendix2} lists the learning, direction,
restart, stagnation, and fallback settings fixed across all reported
experiments. In particular, \(\kappa_d=10^5\) is the direction-scale
setting used in the present implementation and is not claimed to be an
optimal value. The uniform direction safeguard and verified Armijo fallback apply only
to the neural variants. Emergency continuation is recorded separately and
is not an Armijo-certified step.

\subsection{Performance Measures and Profiles}
\label{subsec:performance-measures}

The symbol \(\mathcal P\) denotes the test set, and the common solver set
is defined by
\[
\mathcal S
:=
\left\{
\texttt{HS},
\texttt{PR},
\texttt{HZ},
\texttt{DY},
\texttt{FR},
\texttt{DL},
\texttt{HSNN},
\texttt{PRNN},
\texttt{HZNN},
\texttt{DYNN},
\texttt{FRNN},
\texttt{DLNN}
\right\}.
\]
The same solver set \(\mathcal S\) is used for every performance profile. For consistency with the solver and report labels, this section uses \texttt{PR} and \texttt{PRNN} as abbreviations for \texttt{PRP} and \texttt{PRPNN}, respectively.

For each measure
$m\in\{{\tt nf2g},{\tt ng},{\tt nf},{\tt msec}\}$,
the quantity \(c_{s,p}^{(m)}\) denotes the raw cost of
\(s\in\mathcal S\) on \(p\in\mathcal P\). Lower costs are preferable. Following Dolan--Mor\'e \cite{DolM}, the successful-run ratio is defined by
\begin{equation}
\label{eq:performance-ratio}
r_{s,p}^{(m)}
:=
\frac{
c_{s,p}^{(m)}
}{
\displaystyle
\min_{\bar s\in\mathcal S:\;\bar s\text{ solves }p}
c_{\bar s,p}^{(m)}
}.
\end{equation}
Failed runs have infinite ratios; if no solver succeeds on a problem,
all its ratios are infinite. Successful-run costs are positive.
Then $r_{s,p}^{(m)}\geq1$; equality \(r_{s,p}^{(m)}=1\) means \(s\) is best on
\(p\) for \(m\). For \(\tau\geq1\), the performance profile is defined by
\begin{equation}
\label{eq:performance-profile}
\rho_s^{(m)}(\tau)
:=
\frac{1}{|\mathcal P|}
\left|
\left\{
p\in\mathcal P:
r_{s,p}^{(m)}\leq\tau
\right\}
\right|.
\end{equation}
Thus, \(\rho_s^{(m)}(\tau)\) counts the fraction where \(s\) is within factor \(\tau\)
of the best cost. In particular, $\rho_s^{(m)}(1)$ counts best results for \(s\) under \(m\).
Larger values at \(\tau=1\) indicate more wins.
For \(\tau>1\), \(\rho_s^{(m)}(\tau)\) measures performance within the corresponding factor. Robustness is the limiting height:
\begin{equation}
\label{eq:profile-robustness}
\rho_s^{(m)}(\infty)
:=
\lim_{\tau\rightarrow\infty}
\rho_s^{(m)}(\tau)
={
N_s^{\mathrm{solved}}
}/{
|\mathcal P|
},
\end{equation}
where \(N_s^{\mathrm{solved}}\) counts successful problems for \(s\).
Higher profiles near \(\tau=1\) indicate more best-cost results;
higher limits indicate greater robustness. The comparison uses performance profiles for \({\tt nf2g}\), \({\tt ng}\),
\({\tt nf}\), and \({\tt msec}\) over the relevant \(\tau\) range. No separate
mean-efficiency definition is used.

\subsection{Numerical Results and Pairwise Comparisons}
\label{subsec:summary-comparison}

The six classical/learning-assisted pairs were evaluated on the same
480 unconstrained problems. The numerical data were consolidated from the
individual \texttt{HS}, \texttt{PR}, \texttt{HZ}, \texttt{DY}, \texttt{FR}, and
\texttt{DL} comparison reports. Across the complete solver set \(\mathcal S\),
434 problems were solved by at least one method, while 46 were unsolved by
all 12 methods. Each performance profile uses the common denominator
\(\lvert\mathcal P\rvert=480\) and the minimum cost among successful runs
of all 12 solvers, as prescribed by
\eqref{eq:performance-ratio}--\eqref{eq:performance-profile}.
In particular, problems unsolved by every method remain in the denominator;
they are not discarded as in the separate pairwise source reports.

The solved-count difference is
\begin{equation}
\label{eq:solved-pairwise-change}
\Delta N_{\mathrm{solved}}
:=N_{\texttt{NN}}^{\mathrm{solved}}-N_{\texttt{CG}}^{\mathrm{solved}}.
\end{equation}

\begin{table}[htbp]
\centering
\caption{Successful problems for each classical/learning-assisted pair.
The columns ``\texttt{NN} only'' and ``\texttt{CG} only'' count problems solved exclusively
by the indicated member of the pair. All counts refer to 480 problems.}
\label{tab:pairwise-solved-counts}
\small
\begin{tabular}{llllll}
\toprule
Pair & \texttt{CG} solved & \texttt{NN} solved & \texttt{NN} only & \texttt{CG} only & $\Delta N_{\mathrm{solved}}$\\
\midrule
\texttt{HSNN}-\texttt{HS} & 384 & 414 & 31 & 1 & $+30$\\
\texttt{PRNN}-\texttt{PR} & 397 & 407 & 17 & 7 & $+10$\\
\texttt{HZNN}-\texttt{HZ} & 400 & 413 & 17 & 4 & $+13$\\
\texttt{DYNN}-\texttt{DY} & 357 & 371 & 19 & 5 & $+14$\\
\texttt{FRNN}-\texttt{FR} & 355 & 364 & 18 & 9 & $+9$\\
\texttt{DLNN}-\texttt{DL} & 394 & 412 & 23 & 5 & $+18$\\
\bottomrule
\end{tabular}
\end{table}

Table~\ref{tab:pairwise-solved-counts} shows that every learning-assisted
variant solves more problems than its classical counterpart. The largest
absolute improvement occurs for \texttt{HSNN}, with 30 additional successful
runs; \texttt{DLNN} gains 18, and the remaining gains range from 9 to 14.
The largest solved count is 414 for \texttt{HSNN}, followed by 413 for
\texttt{HZNN} and 412 for \texttt{DLNN}. These improvements are net gains:
there are also classical-only successes in every pair. As is well known, all three {\tt CG} algorithms, namely {\tt HS}, {\tt DL}, and {\tt HZ}, have been developed based on a successive conjugacy condition, independently of the line-search strategy and convexity assumptions. Therefore, the numerical results primarily highlight the effectiveness of the proposed learning-assisted strategy in improving the performance of these methods, which also ensure local conjugacy.

For every measure \(m\), the change in the best-cost fraction and its
percentage-point representation are
\begin{equation}
\label{eq:profile-pairwise-change}
\Delta\rho_m(1):=\rho_{\texttt{NN}}^{(m)}(1)
-\rho_{\texttt{CG}}^{(m)}(1)
\end{equation}
and
\begin{equation}
\label{eq:profile-pairwise-change-pp}
\Delta\rho_m^{\mathrm{pp}}(1):=100\,\Delta\rho_m(1).
\end{equation}
Equal-cost ties are counted for every solver attaining the minimum;
therefore, the best-cost fractions need not sum to one across the methods.

\begin{table}[htbp]
\centering
\caption{Pairwise change of each learning-assisted method relative to its
classical counterpart. The values \(\Delta\rho_m^{\mathrm{pp}}(1)\) are
percentage-point (pp) differences computed with the common 12-solver
reference set and all 480 problems. A positive entry means more best-cost
results, including ties, for the neural method.}
\label{tab:classical-nn-change}
\small
\begin{tabular}{llllll}
\toprule
Pair & $\Delta N_{\mathrm{solved}}$
& $\Delta\rho_{\tt nf2g}^{\mathrm{pp}}(1)$
& $\Delta\rho_{\tt ng}^{\mathrm{pp}}(1)$
& $\Delta\rho_{\tt nf}^{\mathrm{pp}}(1)$
& $\Delta\rho_{\tt msec}^{\mathrm{pp}}(1)$\\
\midrule
\texttt{HSNN}-\texttt{HS} & $+30$ & $+1.04$ & $+1.25$ & $+1.04$ & $-18.33$\\
\texttt{PRNN}-\texttt{PR} & $+10$ & $0.00$ & $0.00$ & $0.00$ & $-23.75$\\
\texttt{HZNN}-\texttt{HZ} & $+13$ & $-0.21$ & $-0.21$ & $-0.42$ & $-20.62$\\
\texttt{DYNN}-\texttt{DY} & $+14$ & $+0.21$ & $+0.62$ & $0.00$ & $-12.50$\\
\texttt{FRNN}-\texttt{FR} & $+9$ & $-0.21$ & $-0.21$ & $-0.21$ & $-15.21$\\
\texttt{DLNN}-\texttt{DL} & $+18$ & $+1.67$ & $+1.88$ & $+1.67$ & $-13.54$\\
\bottomrule
\end{tabular}
\end{table}

\begin{figure}[htbp]
\centering
\begin{subfigure}{.49\textwidth}
\includegraphics[width=\linewidth]{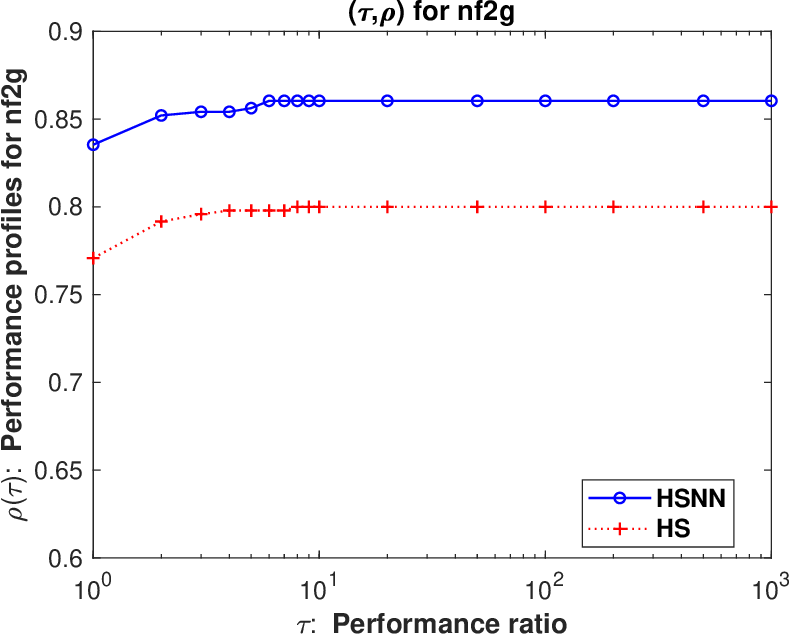}
\caption{\texttt{HS} versus \texttt{HSNN}.}
\label{fig:hs-nf2g}
\end{subfigure}\hfill
\begin{subfigure}{.49\textwidth}
\includegraphics[width=\linewidth]{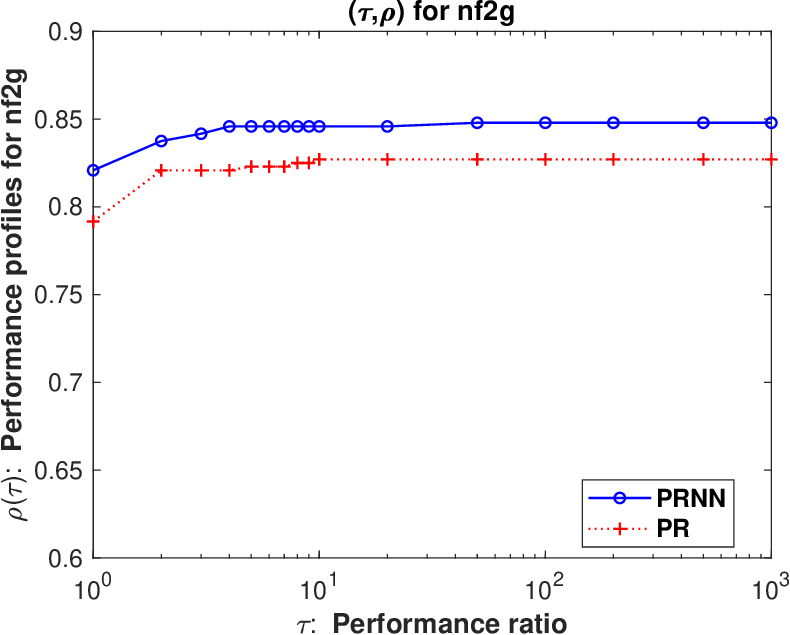}
\caption{\texttt{PR} versus \texttt{PRNN}.}
\label{fig:pr-nf2g}
\end{subfigure}

\begin{subfigure}{.49\textwidth}
\includegraphics[width=\linewidth]{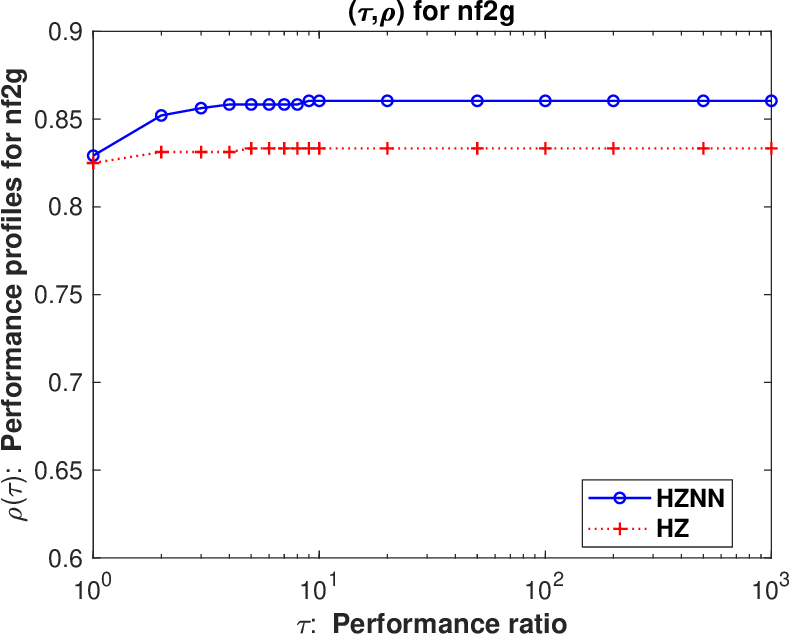}
\caption{\texttt{HZ} versus \texttt{HZNN}.}
\label{fig:hz-nf2g}
\end{subfigure}\hfill
\begin{subfigure}{.49\textwidth}
\includegraphics[width=\linewidth]{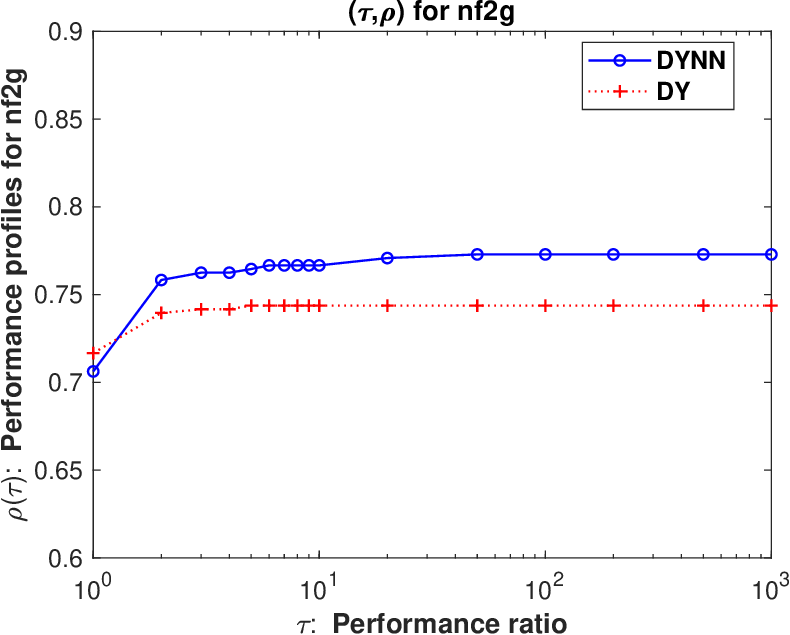}
\caption{\texttt{DY} versus \texttt{DYNN}.}
\label{fig:dy-nf2g}
\end{subfigure}

\begin{subfigure}{.49\textwidth}
\includegraphics[width=\linewidth]{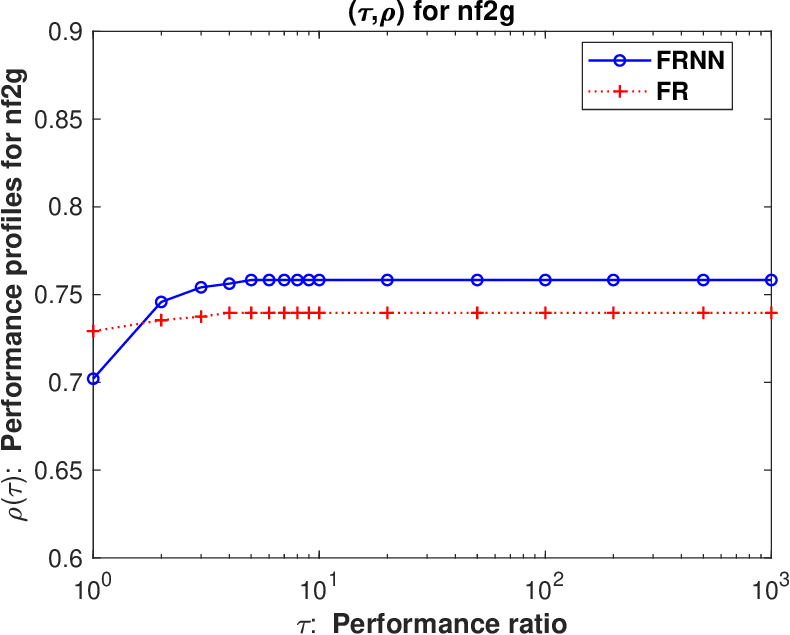}
\caption{\texttt{FR} versus \texttt{FRNN}.}
\label{fig:fr-nf2g}
\end{subfigure}\hfill
\begin{subfigure}{.49\textwidth}
\includegraphics[width=\linewidth]{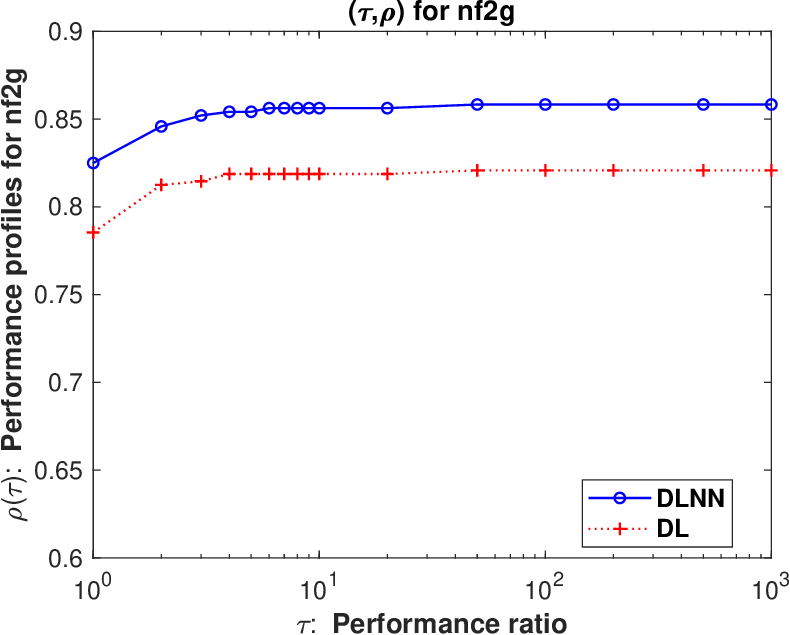}
\caption{\texttt{DL} versus \texttt{DLNN}.}
\label{fig:dl-nf2g}
\end{subfigure}

\caption{Pairwise Dolan--Mor\'e performance profiles for
\texttt{nf2g} on 480 problems, using the common
12-solver baseline.}
\label{fig:pairwise-nf2g-profiles}
\end{figure}

Figure~\ref{fig:pairwise-nf2g-profiles} shows the pairwise \texttt{nf2g}
profiles. Table~\ref{tab:classical-nn-change} reveals modest best-cost
changes despite larger solved-count gains. In particular, \texttt{DLNN}
and \texttt{HSNN} improve their weighted-evaluation win fractions by
1.67 and 1.04 percentage points, respectively. The \texttt{HZNN} and
\texttt{FRNN} variants have slightly smaller best-cost fractions than
their classical counterparts. The profiles therefore do not support a
claim of uniformly reduced evaluation cost.

To separate successful termination from efficiency on problems solved by
both members of a pair, Table~\ref{tab:jointly-solved-comparison} reports
direct cost comparisons. Equal reported costs are listed as ties.

\begin{table}[htbp]
\centering
\caption{Head-to-head results restricted to problems solved by both members
of each pair. Each triple gives the number of problems for which the \texttt{NN}
variant costs less, costs more, or ties its classical counterpart, respectively.}
\label{tab:jointly-solved-comparison}
\small
\begin{tabular}{llll ll l}
\toprule
& \multicolumn{3}{c}{Weighted evaluations \texttt{nf2g}}
& \multicolumn{3}{c}{Elapsed time \texttt{msec}}\\
\cmidrule(lr){2-4}\cmidrule(l){5-7}
Pair & \texttt{NN} less & \texttt{NN} more & Tie & \texttt{NN} less & \texttt{NN} more & Tie\\
\midrule
\texttt{HSNN}-\texttt{HS} & 14 & 13 & 356 & 35 & 307 & 41\\
\texttt{PRNN}-\texttt{PR} & 17 & 13 & 360 & 32 & 315 & 43\\
\texttt{HZNN}-\texttt{HZ} & 4 & 15 & 377 & 22 & 334 & 40\\
\texttt{DYNN}-\texttt{DY} & 13 & 32 & 307 & 30 & 287 & 35\\
\texttt{FRNN}-\texttt{FR} & 5 & 27 & 314 & 34 & 281 & 31\\
\texttt{DLNN}-\texttt{DL} & 17 & 16 & 356 & 43 & 311 & 35\\
\bottomrule
\end{tabular}
\end{table}

Most jointly solved runs have identical evaluation counts, whereas the
classical implementation is usually faster. This difference is visible in
the \texttt{msec} profiles, where every neural variant has a negative
\(\Delta\rho_{\tt msec}^{\mathrm{pp}}(1)\). Measured running times are
sensitive to implementation overhead and their reported resolution; they
should not be interpreted as machine-independent complexity constants.

The six reports provide costs and success indicators but no aggregate
counts of \texttt{ELM} activations, direction repairs, or emergency heuristic
continuations. The comparisons assess the combined learning, safeguarding,
and fallback mechanisms without isolating the contribution of neural
prediction. Runs involving heuristic continuation fall outside the
exact-arithmetic guarantees. The reported results show improved robustness,
with identical weighted evaluation counts on most jointly solved problems.
Storing training data, updating the \texttt{ELM}s, and generating predictions
add computational overhead. The framework is particularly promising
for difficult problems with expensive function or gradient evaluations,
where reductions in evaluation counts, when achieved, may outweigh this
overhead and reduce total running time. To examine this possibility,
Subsection~\ref{subsec:hsnn-conductivity-nf2g} considers a
PDE-constrained conductivity-identification problem and compares
\texttt{HS} with \texttt{HSNN} in terms of the weighted evaluation
measure \texttt{nf2g} and computational time. This example illustrates
a case in which \texttt{HSNN} reduces both evaluation cost and running
time relative to \texttt{HS}.

\subsection{Conductivity Identification with \texttt{NN}-Assisted Step Correction}
\label{subsec:hsnn-conductivity-nf2g}

We consider an inverse-conductivity problem motivated by stationary heat
conduction and electrical conductivity imaging. The goal is to identify
spatial variations of the conductivity from measurements collected under
several source terms. The objective and its exact discrete gradient require
repeated solutions of elliptic partial differential equations. This
example examines whether the safeguards and online predictions of
\texttt{HSNN} can reduce the evaluation cost and elapsed time relative to the classical \texttt{HS}. Synthetic measurements
allow us to check the discrete gradient against finite differences.

Let $\Omega=(0,1)^2$. For each source $q_\ell$, the state $u_\ell$
satisfies
\begin{equation}
\label{eq:cond-nf2g-pde}
-\nabla\cdot\bigl(k(\x;\xi)\nabla u_\ell(\xi)\bigr)
=q_\ell(\xi),\qquad \xi\in\Omega,\qquad
u_\ell|_{\partial\Omega}=0,\qquad \ell=1,\ldots,L.
\end{equation}
Here $\xi$ is a spatial coordinate, whereas $\x\in\R^n$ is the
unconstrained optimization variable. The conductivity is parameterized by
\begin{equation}
\label{eq:cond-nf2g-param}
k(\x;\xi)
=\exp\!\left(\sum_{j=1}^n x_j\phi_j(\xi)\right),
\end{equation}
where the $\phi_j$ are prescribed, spatially localized basis functions.
The exponential representation guarantees positive conductivity without
imposing bounds on $\x$.

Let $N$ denote the number of discrete state unknowns, which need not equal
the optimization dimension $n$. We write $\mathbf{m}=T\x$, with
$T\in\R^{N\times n}$, and obtain the discrete state equations
\begin{equation}
\label{eq:cond-nf2g-discrete}
A(\x)\bu_\ell=\bb_\ell,\qquad
A(\x)=B^{\mathsf T}\operatorname{diag}
\bigl(\exp(RT\x)\bigr)B,\quad \ell=1,\ldots,L.
\end{equation}
The matrix $B$ is the discrete spatial-gradient operator, $R$ maps
nodal log-conductivities to edges, and $\bb_\ell$ is the discrete source
vector. Let $C$ denote the observation operator, $\mathbf{d}_\ell$ the
synthetic observations, and $D$ a spatial-difference matrix. With the
positive weights $s$, $\lambda$, and $\mu$, the inverse problem has the
unconstrained form of \eqref{eq:problem} with
\begin{equation}
\label{eq:cond-nf2g-objective}
f(\x)
:=\frac{s^2}{2}\sum_{\ell=1}^{L}
\norm{C\bu_\ell(\x)-\mathbf{d}_\ell}^{2}
+\frac{\lambda}{2}\norm{DT\x}^{2}
+\frac{\mu}{2}\norm{T\x}^{2}.
\end{equation}

We compute $\g(\x)$ using the exact adjoint of the
\emph{discretized} problem, rather than using a finite-difference
gradient in the optimization. For each load case, the adjoint vector
$\bv_\ell$ solves
\begin{equation}
\label{eq:cond-nf2g-adjoint}
A(\x)\bv_\ell
=s^2C^{\mathsf T}\bigl(C\bu_\ell(\x)-\mathbf{d}_\ell\bigr).
\end{equation}
Define the positive edge conductivity vector
$\mathbf{k}:=\exp(RT\x)$, and let $\odot$ denote componentwise
multiplication. Then, direct differentiation gives
\begin{equation}
\label{eq:cond-nf2g-gradient}
\g(\x)=T^{\mathsf T}\!\left[
-R^{\mathsf T}\!\left(
\mathbf{k}\odot\sum_{\ell=1}^{L}
(B\bu_\ell)\odot(B\bv_\ell)
\right)
+\lambda D^{\mathsf T}D\mathbf{m}+\mu\mathbf{m}
\right].
\end{equation}
Consequently, a joint objective-gradient evaluation assembles and
factorizes a sparse elliptic operator and solves $L$ forward and $L$
adjoint linear systems. A centered directional finite-difference check,
used only for validation, produced a relative error of
$5.88\times10^{-10}$ for the selected instance.

We use a $12\times12$ interior spatial grid ($N=144$ discrete state
unknowns), $n=4$ conductivity parameters, $L=4$ source patterns, and
$100$ scalar observations. The data-misfit scale is $s=20$, with
$\lambda=10^{-10}$ and $\mu=10^{-12}$. Both methods begin from
\[
\x_0=(-2.46746053,\,0.54188715,\,-0.64497769,\,-0.59248754)^{\mathsf T}.
\]
We compare only \texttt{HS} and \texttt{HSNN}, and retain their
original parameter values and line-search settings from
Tables~\ref{tab:cg-wolfe-parameters} and~\ref{tab:nn-parameters}.
For this application experiment, both runs terminate at the stronger
Euclidean criterion $\norm{\g_k}\leq10^{-6}$, which also satisfies
the \texttt{CUTEst} success criterion $\norm{\g_k}_\infty\leq10^{-6}$.
Both runs finish within the common
$\texttt{nf2g}\leq500000$ budget, without a finite time limit.

We use the same evaluation measures as in
Subsection~\ref{subsec:performance-measures}: \texttt{nf2g} and \texttt{msec}. This application driver requests the two
quantities jointly at each trial point, so
$\texttt{nf}=\texttt{ng}$ and $\texttt{nf2g}=3\,\texttt{nf}$.
The elapsed-time measure \texttt{msec} is recorded in milliseconds;
it includes the PDE solves, the line search, and the online
\texttt{ELM} calculations. Reported times are medians of five runs
with alternating solver order.

\begin{figure}[H]
\centering
\includegraphics[width=\linewidth]{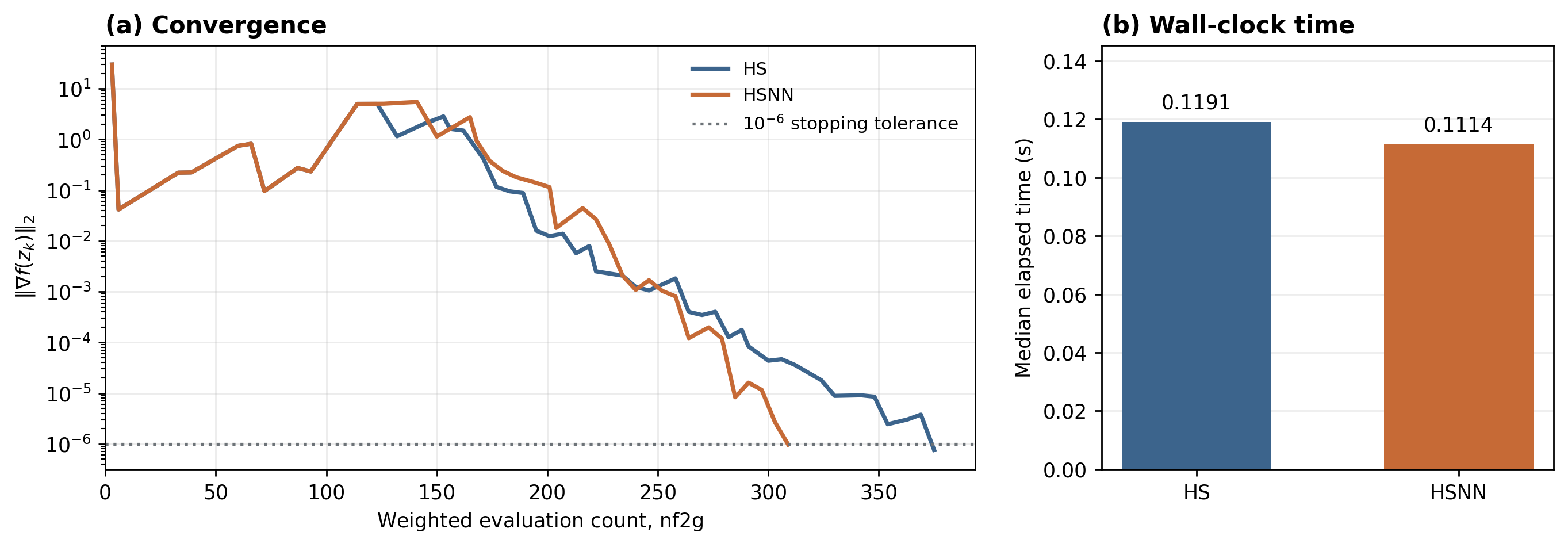}
\caption{Comparison of \texttt{HS} and \texttt{HSNN} on a selected
synthetic inverse-conductivity instance.
\emph{Left:} Euclidean gradient norm against the weighted cost
\texttt{nf2g}; the dotted line is the common $10^{-6}$ tolerance.
\emph{Right:} median elapsed time over five runs, shown in seconds
($0.1191$ and $0.1114$~s, or $119.1$ and $111.4$~\texttt{msec},
respectively). The time measurements include neural-model overhead.}
\label{fig:cond-nf2g-time}
\end{figure}

Figure~\ref{fig:cond-nf2g-time} and
Table~\ref{tab:cond-nf2g-results} show that both methods reach the
required gradient tolerance. Compared with \texttt{HS},
\texttt{HSNN} reduces \texttt{nf2g} from $375$ to $309$, a
$17.6\%$ reduction, and reduces the median measured time from
$119.1$ to $111.4$~\texttt{msec}, approximately $6.5\%$.
Although the evaluation advantage is clear for this instance, the
short elapsed times make the timing difference sensitive to system
noise and implementation details.

\begin{table}[H]
\centering
\caption{Results for the selected inverse-conductivity instance.}
\label{tab:cond-nf2g-results}
\small
\begin{tabular}{llllll}
\toprule
Method & Iterations & \texttt{nf} & \texttt{ng} &
\texttt{nf2g} & \texttt{msec} \\
\midrule
\texttt{HS}   & 46 & 125 & 125 & 375 & 119.1 \\
\texttt{HSNN} & 36 & 103 & 103 & \textbf{309} & \textbf{111.4} \\
\bottomrule
\end{tabular}
\end{table}

The recorded solver trace identifies an intervention by the
fallback-step \texttt{ELM} described in
the line-search-failure mechanism in Section~\ref{sec:nn-parameters}. After an unsuccessful
\texttt{cvsrch} call, the trained model proposes enlarging the initial
fallback step from $2.72408\times10^{-6}$ to
$2.72408\times10^{-2}$. The accepted step changes the subsequent
optimization trajectory. The intervention is therefore a learned
\emph{step correction} following a line-search failure, not a
learned correction of the \texttt{HS} coefficient. No trained
coefficient-\texttt{ELM} direction was accepted during this run.

The Hessian condition number at the reference parameter for this
selected instance is approximately $5.27\times10^1$. Hence this test
isolates a helpful learned fallback but does not demonstrate neural
recovery on a severely ill-conditioned parameter-identification
problem. A separate experiment with condition number about
$1.08\times10^6$ gave \texttt{nf2g} costs of $1491$ for \texttt{HS}
and $543$ for \texttt{HSNN}, but its advantage came from deterministic
bootstrap recovery, not an accepted neural prediction. Additional
completed $8$- and $16$-parameter examples had equal weighted
evaluation counts because the learned corrections were not used.
The favorable four-parameter instance was identified through
exploratory screening. Its evaluation and time improvements should
therefore be viewed as a proof of mechanism rather than evidence of
uniform or large-scale superiority. A prespecified application test
set and a native \texttt{MATLAB} replication would be needed for a broader
performance claim.

\section{Concluding Remarks}
\label{sec:conclusion}

This study has developed a safeguarded learning-assisted framework for
nonlinear \texttt{CG} methods in which online prediction supports, rather
than replaces, the underlying classical optimization process. The standard
two-term \texttt{CG} structure is retained, while two independent online
\texttt{ELM}s are used only at selected intervention points. The first may
provide a coefficient during restart recovery or persistent stagnation, and
the second may modify the initial fallback-step proposal after a line-search
failure. Deterministic descent, scale, distributional, and fallback
safeguards govern the use of these predictions.

The numerical experiments on 480 \texttt{CUTEst} problems show improved
robustness for all six learning-assisted variants. Each learning-assisted
method solves more problems than its corresponding classical method, with
the gains ranging from 9 to 30 additional solved problems. The largest
increase is obtained by \texttt{HSNN}, which solves 414 problems compared
with 384 for \texttt{HS}. Across the test set, the evaluation costs remain
broadly comparable between the learning-assisted and classical variants.

The additional robustness is not computationally free. Feature
construction, archive maintenance, \texttt{ELM} fitting, prediction, and
the deterministic safeguards introduce algebraic overhead beyond that of
the corresponding classical \texttt{CG} method. Therefore, the results do
not support a general claim that the learning-assisted variants reduce
wall-clock time. Their potential computational advantage is instead most
relevant when function or gradient evaluations are expensive. In such
settings, reductions in evaluation counts, when achieved, may outweigh the
learning overhead and consequently reduce total computational time.

The PDE-constrained conductivity-identification experiment illustrates this
possibility. On the selected instance, the useful learned intervention is a
fallback-step correction after a line-search failure rather than a learned
\texttt{CG} coefficient. The observed reduction in evaluation cost and
runtime should therefore be interpreted as evidence that the proposed
mechanism can be beneficial when evaluations are expensive, rather than as
evidence of uniform superiority over the classical methods.

A central feature of the framework is that predictions remain subordinate
to deterministic optimization safeguards. An unavailable,
out-of-distribution, nonfinite, or inadmissible prediction does not by
itself determine the next iterate; instead, the appropriate deterministic
recovery or fallback mechanism remains available. The verified core of the
method is covered by the convergence and evaluation-complexity analysis,
whereas the optional emergency continuation is explicitly excluded from
those guarantees.

Finally, the observed gains arise from the combined action of neural
prediction, deterministic safeguarding, restart recovery, stagnation
handling, and line-search fallback mechanisms. The present experiments do
not isolate the individual contribution of each component. Dedicated
ablation studies are therefore an important direction for future work,
together with tests on broader and more expensive problem classes and
further investigation of learning-assisted strategies that preserve the
low-storage character of nonlinear \texttt{CG} methods.

{\bf Funding.} Morteza Kimiaei acknowledges financial support from the Austrian Science Fund under \url{https://doi.org/10.55776/PAT2747625}.

%{\bf Declaration of interests.} The authors declare that they have no competing interests.

%{\bf CRediT author statement.} Morteza Kimiaei and Saman Babaie--Kafaki contributed equally to the Conceptualization, Methodology, Formal analysis, Investigation, Validation, Writing -- original draft, and Writing -- review \& editing of this work.

%%%%%
%%%%%
%%%%%
\appendix

\section{Classical Method and Line-Search Parameters}\label{sec:appendix1}

\footnotesize
\begin{longtable}{lll}
\caption{Classical \texttt{CG}, stopping, and Wolfe line-search parameters
prescribed for the experiments.}
\label{tab:cg-wolfe-parameters}\\
\toprule
Setting & Value & Role \\
\midrule
\endfirsthead

\toprule
Setting & Value & Role \\
\midrule
\endhead

\bottomrule
\endfoot

\bottomrule
\endlastfoot

\(t_{\mathrm{DL}}\) & \(0.1\) &
fixed \texttt{DL} parameter for \texttt{DL}/\texttt{DLNN} \\
\texttt{DY}/\texttt{FR} restart threshold & \(\sqrt{\epsilon_{\mathrm{mach}}}\) &
criterion in \eqref{eq:dy-fr-restart} \\
\(\delta_a\) & \(10^{-8}\) &
original fallback-step seed factor \\
Stopping tolerance & \(10^{-6}\) &
\(\norm{\g}_\infty\) tolerance \\
Evaluation budget & \(500000\) &
maximum \({\tt nf2g}\) \\
Time limit & \(\infty\) &
no time-based termination \\
\midrule
Wolfe initial step & \(1\) &
prescribed \texttt{cvsrch} initial step \\
\(c_1\) & \(10^{-4}\) &
sufficient-decrease parameter \\
\(c_2\) & \(0.1\) &
Wolfe curvature parameter \\
\(c_4\) & \(10^{-3}\) &
declared; inactive in the acceptance test \\
Maximum line-search evaluations \(\ell_{\max}\) & \(30\) &
solver-enforced evaluation limit \\
Minimum line-search step & \(0\) &
prescribed \texttt{cvsrch} lower bound \\
Maximum line-search step & \(\infty\) &
prescribed \texttt{cvsrch} upper bound \\
Bracket tolerance & \(10^{-9}\) &
prescribed \texttt{cvsrch} tolerance \\
\end{longtable}

\normalsize

\section{Learning-Assisted Method Parameters}\label{sec:appendix2}

\scriptsize
\begin{longtable}{lll}
\caption{Parameters of the learning-assisted \texttt{CGNN} variants.}
\label{tab:nn-parameters}\\
\toprule
Parameter & Value & Role \\
\midrule
\endfirsthead

\toprule
Parameter & Value & Role \\
\midrule
\endhead

\bottomrule
\endfoot

\bottomrule
\endlastfoot

\multicolumn{3}{l}{{\bf Coefficient-\texttt{ELM} and recovery}}\\
Input dimension \(p_\beta\) & \(12\) &
coefficient feature dimension \\
Hidden units \(h_\beta\) & \(24\) &
coefficient-\texttt{ELM} hidden layer \\
Archive size \(N_\beta^{\max}\) & \(250\) &
maximum replay samples \\
Minimum training samples & \(8\) &
first eligible fit \\
Retraining interval & \(4\) &
new samples between fits \\
Ridge parameter & \(10^{-4}\) &
regularized least squares \\
Random seed & \(73\) &
fixed hidden-layer initialization \\
Standard-deviation threshold \(\sigma_{\min}\) & \(10^{-12}\) &
standardization safeguard threshold \\
Feature clipping \(c_{\mathrm{feat}}\) & \(8\) &
feature clipping bound \\
\texttt{OOD} factor \(\gamma_{\texttt{OOD}}\) & \(5\) &
\texttt{OOD} threshold factor \\
\texttt{OOD} archive quantile \(\nu_{\texttt{OOD}}\) & \(0.95\) &
quantile used to construct the \texttt{OOD} archive scale \\
Coefficient history length & \(12\) &
history-bootstrap memory \\
\(\beta_{\max}\) & \(10\) &
neural/recovery coefficient cap \\
\(\varepsilon_\beta\) & \(10^{-8}\) &
coefficient-target clipping away from \(\pm1\) \\
Zero-coefficient tolerance & \(0\) &
zero-coefficient trigger \\
\(\beta_{\min}^{\mathrm{nz}}\) & \(10^{-10}\) &
nonzero-recovery threshold \\
Coefficient-repair contraction \(\vartheta_\beta\) & \(0.5\) &
contraction used by local descent repair \\
Maximum repair contractions & \(24\) &
finite-precision recovery limit \\

\midrule
\multicolumn{3}{l}{{\bf Uniform search-direction and restart safeguards}}\\
\(\Delta_{\texttt{SDC}}=\delta\) & \(\epsilon_{\mathrm{mach}}\) &
local/final sufficient-descent setting \\
\(\kappa_d\) & \(10^5\) &
current direction-scale bound \\
Finite-precision DSC factor \(\chi_{\mathrm{fp}}\) & \(100\) &
inflation factor multiplying \(\epsilon_{\mathrm{mach}}\) \\
\(\rho\) & \(1\) &
\texttt{DYNN}/\texttt{FRNN} inherited-component bound \\

\midrule
\multicolumn{3}{l}{{\bf Persistent-stagnation mechanism}}\\
History window & \(5\) &
objective/gradient stagnation window \\
\(\tau_f\) & \(10^{-3}\) &
objective stagnation tolerance \\
\(\tau_g\) & \(10^{-3}\) &
gradient stagnation tolerance \\
First eligible loop iteration & \(20\) &
one-based intervention index \\
Spacing between attempts & \(5\) &
minimum iteration separation \\
\(\eta_{\mathrm{dir}}\) & \(0.005\) &
minimum tentative direction change \\

\midrule
\multicolumn{3}{l}{{\bf Fallback-step \texttt{ELM} and Armijo fallback safeguard}}\\
Input dimension \(p_\alpha\) & \(10\) &
fallback-step feature dimension \\
Hidden units \(h_\alpha\) & \(16\) &
fallback-step hidden layer \\
Archive size \(N_\alpha^{\max}\) & \(250\) &
maximum replay samples \\
Minimum training samples & \(8\) &
first eligible fit \\
Retraining interval & \(4\) &
new samples between fits \\
Ridge parameter & \(10^{-4}\) &
regularized least squares \\
Random seed & \(1082\) &
\(73+1009\) \\
Standard-deviation threshold \(\sigma_{\min}\) & \(10^{-12}\) &
standardization safeguard threshold \\
Feature clipping \(c_{\mathrm{feat}}\) & \(8\) &
standardized feature bound \\
\texttt{OOD} factor \(\gamma_{\texttt{OOD}}\) & \(5\) &
\texttt{OOD} threshold factor \\
\texttt{OOD} archive quantile \(\nu_{\texttt{OOD}}\) & \(0.95\) &
quantile used to construct the \texttt{OOD} archive scale \\
\(C_{\min}\) & \(1\) &
minimum neural correction factor \\
\(C_{\max}\) & \(10^4\) &
maximum neural correction factor \\
\(C_{\mathrm{tar}}\) & \(10^8\) &
clipping ratio used for fallback-step training targets \\
\(\alpha_{\min}\) & \(10^{-8}\) &
minimum repaired initial fallback step \\
\(\alpha_{\max}\) & \(1\) &
maximum initial fallback step \\
\(\vartheta\) & \(0.5\) &
Armijo backtracking contraction \\
Emergency continuation & enabled &
practical option outside the verified analysis \\
\(\epsilon_{\mathrm H}\) & \(10^{-8}\) &
scaled objective allowance above the best value \\
\end{longtable}
\end{sloppypar}

\normalsize

%%%%%
%%%%%
%%%%%
%\bibliographystyle{elsarticle-num}
%\bibliography{Refs}

\end{document}